\documentclass[11pt]{article}
\usepackage{amsfonts}
\usepackage{color}
\usepackage{amsmath,amssymb,comment}
\newtheorem{theo}{Theorem}[section]
\newtheorem{lem}[theo]{Lemma}

\newcommand{\mysection}[1]{\section{#1} \setcounter{equation}{0}}
\newcommand{\proof}{{\sc Proof.} \quad}
\newcommand{\proofc}{{\sc Proof} \ }
\newcommand{\be}{\begin{equation} \label}
\newcommand{\ee}{\end{equation}}
\newcommand{\bea}{\begin{eqnarray}\label}
\newcommand{\eea}{\end{eqnarray}}
\newcommand{\bas}{\begin{eqnarray*}}
\newcommand{\eas}{\end{eqnarray*}}
\newcommand{\bit}{\begin{itemize}}
\newcommand{\eit}{\end{itemize}}
\newcommand{\qed}{\hfill$\Box$ \vskip.2cm}
\newcommand{\nn}{\nonumber}
\newcommand{\R}{\mathbb{R}}
\newcommand{\N}{\mathbb{N}}
\newcommand{\pO}{\partial\Omega}

\newcommand{\eps}{\varepsilon}

\newcommand{\hra}{\hookrightarrow}
\newcommand{\io}{\int_\Omega}
\newcommand{\na}{\nabla}
\newcommand{\Del}{\Delta}
\newcommand{\del}{\delta}
\newcommand{\al}{\alpha}
\newcommand{\lam}{\lambda}

\newcommand{\pa}{\partial}
\newcommand{\bom}{\overline{\Omega}}
\newcommand{\Om}{\Omega}

\newcommand{\ov}{\overline}

\newcommand{\wh}{\widehat}

\newcommand{\hs}{\hspace*}

\newcommand{\vp}{\varphi}
\newcommand{\lbal}{\left\{ \begin{array}{l}}
\newcommand{\lball}{\left\{ \begin{array}{ll}}
\newcommand{\ear}{\end{array} \right.}

\newcommand{\abs}{\\[5pt]}
\newcommand{\Abs}{\\[5mm]}

\newcommand{\adb}{\allowdisplaybreaks}
\newcommand{\tm}{T_{max}}

\newcommand{\oy}{\ov{y}}

\newcommand{\ovv}{\ov{v}}

\renewcommand{\div}{{\rm div} \,}
\begin{document}
\adb
\title{
Describing smooth small-data solutions to a\\
quasilinear hyperbolic-parabolic system\\
by $W^{1,p}$ energy analysis}
\author{
Leander Claes\footnote{claes@emt.uni-paderborn.de}\\
{\small Universit\"at Paderborn}\\
{\small Institut f\"ur Elektrotechnik und Informationstechnik}\\
{\small 33098 Paderborn, Germany}
\and
Michael Winkler\footnote{michael.winkler@math.uni-paderborn.de}\\
{\small Universit\"at Paderborn}\\
{\small Institut f\"ur Mathematik}\\
{\small 33098 Paderborn, Germany}}
\date{}
\maketitle
\begin{abstract}
\noindent
In bounded $n$-dimensonal domains with $n\ge 1$, this manuscript examines an initial-boundary value problem for
the system
\bas
	\lbal
	u_{tt} = \nabla \cdot (\gamma(\Theta) \na u_t) + a \na \cdot (\gamma(\Theta) \na u) + \na\cdot f(\Theta), \\[1mm]
	\Theta_t = D\Del\Theta + \Gamma(\Theta) |\na u_t|^2 + F(\Theta)\cdot \na u_t,
	\ear
\eas
which in the case $n=1$ and with $\gamma\equiv \Gamma$ as well as $f\equiv F$ reduces to
the classical model for the evolution of strains and temperatures in thermoviscoelasticity.
Unlike in previous related studies, the focus here is on situations in which besides $f$ and $F$, also the core ingredients
$\gamma$ and $\Gamma$ may depend on the temperature variable $\Theta$.\abs
Firstly, a statement on local existence of classical
solutions is derived for arbitrary $a>0, D>0$ as well as $0<\gamma\in C^2([0,\infty))$ and $0\le\Gamma\in C^1([0,\infty))$,
for functions
$f\in C^2([0,\infty);\R^n)$ and
$F\in C^1([0,\infty);\R^n)$ with $F(0)=0$,
and for suitably regular initial data of arbitrary size.
Secondly, it is seen that for each $p\ge 2$ such that $p>n$ there exists $\del(p)>0$ with the property that whenever
in addition to the above we have
\bas
	\frac{a}{\gamma(0)} \le \del(p)
	\qquad \mbox{and} \qquad
	\frac{|f'(\Theta_\star)| \cdot |F(\Theta_\star)|}{D \cdot \gamma(\Theta_\star)}
	\le \del(p),
\eas
for initial data suitably close to the constant level given by $u=0$ and $\Theta=\Theta_\star$, with any fixed
$\Theta_\star\ge 0$,
these solutions are actually global in time and have the property that $\na u_t, \na u$ and $\na\Theta$
decay exponentially fast in $L^p$.
This is achieved by detecting suitable dissipative properties of functionals involving norms
of these gradients in $L^p$ spaces.\abs
\noindent {\bf Key words:} nonlinear acoustics; thermoviscoelasticity; small-data solution; decay rate\\
{\bf MSC 2020:} 35B40 (primary); 35L05, 35B35 (secondary)
\end{abstract}
\newpage
\section{Introduction}\label{intro}
In a smoothly bounded domain $\Om\subset\R^n$ with $n\ge 1$, this manuscript is concerned with the initial-boundary
value problem
\be{0}
	\lball
	u_{tt} = \nabla \cdot (\gamma(\Theta) \na u_t) + a \na \cdot (\gamma(\Theta) \na u) + \na\cdot f(\Theta),
	\qquad & x\in\Om, \ t>0, \\[1mm]
	\Theta_t = D\Del\Theta + \Gamma(\Theta) |\na u_t|^2 + F(\Theta)\cdot\na u_t,
	\qquad & x\in\Om, \ t>0, \\[1mm]
	\frac{\pa u}{\pa\nu}=\frac{\pa\Theta}{\pa\nu}=0,
	\qquad & x\in\pO, \ t>0, \\[1mm]
	u(x,0)=u_0(x), \quad u_t(x,0)=u_{0t}(x), \quad \Theta(x,0)=\Theta_0(x),
	\qquad & x\in\Om,
	\ear
\ee
where $0\le \gamma$ and $0\le \Gamma$, $f$ and $F$ are prescribed functions on $[0,\infty)$,
where $a$ and $D$ are positive parameters, and where
$u_0$, $u_{0t}$ and $\Theta_0$ are given and suitably regular functions satisfying $\Theta_0\ge 0$.
As described in detail in \cite{roubicek},
when reduced by choosing $\gamma\equiv \Gamma$, $f\equiv F$ and $n=1$,
problems of this form arise in the modeling of heat generation by acoustic waves in
viscoelastic materials of one-dimensional structure,
based on the assumption that mechanical losses occur according to a Kelvin-Voigt type material law
(\cite{GutierrezLemini2014}, \cite{Meyers2008}).\abs
In fact, in the case when $n=1$ the PDE system in (\ref{0}) precisely coincides with the corresponding one-dimensional
model thereby obtained, and then $u=u(x,t)$ and $\Theta=\Theta(x,t)$ represent the scalar displacement field and the
temperature distribution inside the material, respectively;
in its higher-dimensional version,
(\ref{0}) may be viewed as a simplified variant of a more complex model which
describes the vector-valued displacement variable $u(x,t)\in\R^n$, along with the associated temperature field,
on the basis of the system
\be{01}
	\lbal
	u_{tt} = d \div (\gamma(\Theta) : \na^s u_t) + a\div (\gamma(\Theta) : \na^s u) + \div f(\Theta), \\[1mm]
	{\mathcal{C}}(\Theta)
	\Theta_t = D\Del\Theta + d \langle \gamma(\Theta) : \na^s u_t , \na^s u_t \rangle + \langle f(\Theta) , \na^s u \rangle,
	\ear
\ee
with the given fourth-rank elasticity tensor $\gamma=\gamma(\Theta)$, and with $\na^s u$ denoting the symmetric gradient of $u$
(\cite{roubicek}, \cite{claes_lankeit_win}).\abs
In its higher-dimensional version, this latter system seem yet far from understood already in the case when
$\gamma\equiv const.$; in fact, despite remarkably deep analysis undertaken for various closely related systems
(\cite{mielke_roubicek}, \cite{rossi_roubicek_interfaces13}, \cite{yoshikawa_pawlow_zajaczkowski_SIMA}, \cite{roubicek_nodea13},
\cite{owczarek_wielgos}, \cite{zimmer}),
results on solvability for versions of (\ref{01}) with $n\ge 2$ either concentrate on local-in-time
(\cite{jiang_racke90}, \cite{bonetti_bonfanti}) or small-data solutions (\cite{racke90}, \cite{shibata}),
or rely on further restrictions such as hypotheses on
suitably mild and especially sublinear growth of $\Theta\mapsto f(\Theta)$
(\cite{blanchard_guibe}, \cite{roubicek}),
on appropriately strong growth of the temperature-dependent heat capacity ${\mathcal{C}}$ with respect to $\Theta$
(\cite{roubicek}, \cite{pawlow_zajaczkowski_cpaa17}, \cite{gawinecki_zajaczkowski_cpaa}, \cite{gawinecki_zajaczkowski}),
or on the absence of viscoelastic contributions,
through requiring that $d=0$ inter alia neglecting the crucial quadratic nonlinear heat production term
$d \langle \gamma(\Theta) : \na^s u_t , \na^s u_t \rangle$ (\cite{cieslak}).\abs
In comparison to these results which mainly resort to suitably weakened concepts of solvability,
corresponding one-dimensional studies go significantly further by constructing classical solutions
in a large variety of particular scenarios (\cite{slemrod}, \cite{dafermos_hsiao_smooth}, \cite{kim},
\cite{jiang1990},\cite{chen_hoffmann},
\cite{chen_hoffmann2}, \cite{racke_zheng_JDE1997}, \cite{song_jiang}, \cite{dafermos}).
Common to these seems that explicit temperature dependencies are considered mainly in the model parameters related to thermoelastic effects;
experimental observations reveal, however, that the core ingredient $\gamma$, which quantifies the elastic behavior
may considerably vary with $\Theta$ (\cite{friesen}, \cite{Gubinyi2007}).
Further, modelling the viscoelastic properties using a scalar factor for the dissipative term in the viscous wave equation significantly reduces the number of parameters necessary for a quantitative, three-dimensional description and is thus advantageous in for application purposes (\cite{Lemaitre2001}, \cite{Carcione2004}).
The most common damping model used in numerical finite element methods (\cite{Craig2006})
also employs scalar factors to scale the dissipative terms in the differential equations to be solved.
Quantitative methods of material identification show that physical observations can be reproduced using this model
(\cite{Kim2009}, \cite{Feldmann2021}, \cite{Claes2023}).
For processes with higher frequencies, this model coincides with elastic part of the system considered here.
However, the literature concerned with associated versions of (\ref{01}) seems to essentially reduce to a recent result on global generalized solvability in the case when ${\mathcal{C}}\equiv const.$ and $f\equiv 0$ (\cite{claes_lankeit_win}).\abs
{\bf Main results.} \quad
The main purpose of the present manuscript consists in undertaking a basic step into the understanding of (\ref{0})
in the presence of fairly general temperature-dependent constituents $\gamma, \Gamma, f$ and $F$.
From a methodological perspective, our focus will be on presenting an approach
which 	
does not rely on evolution properties of high-order energy functionals, but which is nevertheless
capable of carving out essential features of (\ref{0}) in frameworks of 	
suitably small data.
Since this approach will turn out to be independent of dimensionality in essential part,
in view of a possible potential for extensions to related problems from nonlinear acoustics in liquids and gases
in which an analysis can concentrate on scalar quantities (cf., e.g., \cite{lasiecka_limit_model},
\cite{kaltenbacher_lasiecka_pospieszalska}, \cite{lasiecka_wang}, \cite{lasiecka_wang2}),
\cite{kaltenbacher_nikolic_M3AS2019}, \cite{kaltenbacher_nikolic_SIMA} for some representative studies)
we do not restrict our considerations to the case $n=1$ of immediate physical relevance, but rather include
higher-dimensional situations as well.\abs
In line with our ambition to remain within the realm of smooth solutions throughout,
as a prerequisite for this the following first of our results asserts local classical solvability along with a
handy extensibility criterion, to be achieved in Section \ref{sect2} on the basis of Schauder's fixed point theorem
and parabolic regularity arguments applied to the Neumann problem for the system
\be{equiv}
	\lbal
	v_t = \na\cdot (\gamma(\Theta)\na v) + av - a^2 u + \na \cdot f(\Theta), \\[1mm]
	u_t = v-au, \\[1mm]
	\Theta_t = D\Del\Theta + \Gamma(\Theta) |\na v - a\na u|^2 + F(\Theta)\cdot (\na v-a\na u),
	\ear
\ee
to which (\ref{0}) becomes equivalent upon letting $v:=u_t+au$. Indeed:
\begin{theo}\label{theo_loc}
  Let $n\ge 1$ and $\Om\subset\R^n$ be a bounded domain with smooth boundary, and suppose that
  \be{gf}
	\lbal
	\mbox{$a>0$ and $D>0$ are constants, that} \\[1mm]
	\mbox{$\gamma\in C^2([0,\infty))$ and $\Gamma\in C^1([0,\infty))$ are such that $\gamma>0$ and $\Gamma\ge 0$ on $[0,\infty)$, and that} \\[1mm]
  	\mbox{$f\in C^2([0,\infty);\R^n)$ and $F\in C^1([0,\infty);\R^n)$ satisfy $F(0)=0$.}
	\ear
  \ee
  Then whenever
  \be{Init}
	\lbal
	u_0\in C^2(\bom)
	\mbox{ is such that $\frac{\pa u_0}{\pa\nu}=0$ on $\pO$,} \\[1mm]
	u_{0t}\in C^{1+\al}(\bom)		
	\mbox{ is such that $\frac{\pa u_{0t}}{\pa\nu}=0$ on $\pO$, \qquad and} \\[1mm]
	\Theta_0\in C^{1+\al}(\bom)
	\mbox{ satisfies $\Theta_0\ge 0$ in $\Om$ and $\frac{\pa \Theta_0}{\pa\nu}=0$ on $\pO$,}
	\ear
  \ee
  with some $\al\in (0,1)$,
  there exist $\tm\in (0,\infty]$ as well as functions
  \be{tl1}
	\lbal
	u\in \Big( \bigcup_{\beta\in (0,1)} C^{1+\beta,\frac{1+\beta}{2}}(\bom\times [0,\tm))\Big) \cap C^{2,1}(\bom\times (0,\tm))
		\qquad \mbox{and} \\[1mm]
	\Theta\in \Big( \bigcup_{\beta\in (0,1)} C^{1+\beta,\frac{1+\beta}{2}}(\bom\times [0,\tm))\Big)
		\cap C^{2,1}(\bom\times (0,\tm))
	\ear
  \ee
  which are such that
  \be{tl2}
	\begin{array}{l}
	u_t\in
	\Big( \bigcup_{\beta\in (0,1)} C^{1+\beta,\frac{1+\beta}{2}}(\bom\times [0,\tm))\Big) \cap C^{2,1}(\bom\times (0,\tm)),
	\end{array}
  \ee
  that $\Theta\ge 0$ in $\Om\times (0,\tm)$, that $(u,\Theta)$ solves (\ref{0}) in the classical pointwise sense
  in $\Om\times (0,\tm)$, and which have the additional property that
  \bea{Ext}
	& & \hs{-15mm}
	\mbox{if $\tm<\infty$, \quad then \quad} \nn\\
	& & \hs{-6mm}
	\limsup_{t\nearrow\tm} \Big\{ \|u_t(\cdot,t)\|_{W^{1,p}(\Om)} + \|\Theta(\cdot,t)\|_{L^\infty(\Om)} \Big\}
		=\infty
	\mbox{\quad for all $p\ge 2$ such that $p>n$.}
  \eea
\end{theo}
{\bf Remark.} \quad
i) \ The criterion in (\ref{Ext}) asserts extensibility under rather assumptions on solution regularity which appear relatively
mild in comparison to related precedents achieved by analysis of Hilbert space type energies.
In fact, the local theory accordingly developed in \cite{racke_shibata_zheng} for problems of the form
\bas
	\lbal
	u_{tt} = A_1(u_x,\Theta)u_{xx} + A_2(u_x,\Theta) \Theta_x + f(x,t), \\[1mm]
	A_3(u_x,\Theta) \Theta_t = D(\Theta,\Theta_x)\Theta_{xx} + A_2(u_x,\Theta) u_{xt} + g(x,t),
	\ear
\eas
involves higher differentiability
of the initial data, and hence can be seen to imply that for the solutions correspondingly obtained there,
\bas
	& & \hs{-20mm}
	\mbox{if $\tm<\infty$, \quad then \quad}
	\lim_{t\nearrow\tm} \Big\{ \|u(\cdot,t)\|_{W^{3,2}(\Om)} + \|u_t(\cdot,t)\|_{W^{2,2}(\Om)} + \|u(\cdot,t)\|_{W^{1,2}(\Om)}
		\\
	& & \hs{60mm}
	+ \|\Theta(\cdot,t)\|_{W^{3,2}(\Om)} + \|\Theta_t(\cdot,t)\|_{W^{2,2}(\Om)} \Big\} = \infty.
\eas
Similarly, also in the local results achieved in \cite{nikolic_saidhouari_NON} for the Westervelt-Pennes system
\bas
	\lbal
	u_{tt} = b \Del u_t + a(\Theta) \Del u + k(\Theta) (u^2)_{tt}, \\[1mm]
	\Theta_t = D\Del\Theta - A\Theta + B u_t^2,
	\ear
\eas
from nonlinear thermoacoustics in liquid media require a control of, inter alia, the spatial $H^3$ norm of $u$ for prolongation.\\
ii) \ We emphasize that no assumptions on upper nor uniform lower bounds for $\gamma$ are required in Theorem \ref{theo_loc}.
Accodingly, the above statement in particular also covers singular and degenerate cases obtained by choosing unbounded
functions $\gamma$, or by considering coefficients $\gamma$ for which $\gamma(\xi)>0$ for all $\xi\ge 0$, but which
satisfy $\gamma(\xi)\to 0$ as $\xi\to\infty$.\Abs
Section \ref{sect3} will thereafter contain the core of our analysis, by namely investigating the time evolution of
\be{en1}
	y(t):=\io |\na v|^p + C_1 \io |\na u|^p + C_2 \io |\na\Theta|^p + C_3 \io |\na u|^{p+2}
\ee
for arbitrary $p\ge 2$, fixed temperature background levels $\Theta_\star\ge 0$ and appropriately chosen
$C_i=C_i(p,a,D,\gamma,\Gamma,f,F,\Theta_\star)>0$, $i\in\{1,2.3\}$.
Indeed, relying on a Poincar\'e type functional inequality of the form
\be{poincare}
	\io |\na\vp|^p \le C(p) \io |\na\vp|^{p-2} |D^2\vp|^2,
\ee
valid with some $C(p)>0$ for all $\vp\in C^2(\bom)$ fulfilling $\frac{\pa\vp}{\pa\nu}=0$ on $\pO$ (Lemma \ref{lem6}),
Lemma \ref{lem9} will show that whenever $a$ is sufficiently small, this function satisfies
\be{energy}
	y'(t) + c_1 y(t) + c_1 \cdot \bigg\{ \io |\na v|^{p-2} |D^2 v|^2 + \io |\na \Theta|^{p-2} |D^2\Theta|^2 \bigg\}
	\le c_2 y^\lam(t)
\ee
throughout evolution, where $c_i=c_i(p,a,D,\gamma,\Gamma,f,F,\Theta_\star)>0$,
and where the inequality $\lam=\lam(p)>1$ ensures that $y$
remains below its initial value $y(0)$ provided that the latter is small enough.
Combining this with the information on second-order dissipation, as additionally contained in (\ref{energy}), will
lead to our main result concerning the behavior of (\ref{0}) in small-data settings:
\begin{theo}\label{theo12}
  Let $n\ge 1$ and $\Om\subset\R^n$ be a smoothly bounded domain, and let $p\ge 2$ be such that $p>n$.
  Then there exists $\del(p)>0$ with the property that whenever
  and whenever $a,D\gamma,\Gamma,f$ and $F$ are such that (\ref{gf}) holds and that
  \be{12.1}
	\frac{a}{\gamma(\Theta_\star)} \le \del(p)
  \ee
  as well as
  \be{9.01}
	\frac{|f'(\Theta_\star)| \cdot |F(\Theta_\star)|}{D \cdot \gamma(\Theta_\star)}
	\le \del(p)	
  \ee
  with some $\Theta_\star\ge 0$,
  it is possible to choose $\eps_\star=\eps_\star(p,a,D,\gamma,\Gamma,f,F,\Theta_\star)>0$,
  $\kappa=\kappa(p,a,D,\gamma,\Gamma,f,F,\Theta_\star)>0$,
  $(\eta_\star(\eps))_{\eps\in (0,\eps_\star)} \subset (0,\infty)$ and
  $C=C(p,a,D,\gamma,\Gamma,f,F,\Theta_\star)>0$
  in such a way that if $(u_0,u_{0t},\Theta_0)$ satisfies (\ref{Init})
  with
  \be{12.3}
	\|\Theta_0-\Theta_\star\|_{L^\infty(\Om)} \le \eps
  \ee
  \be{12.2}
	\io |\na u_{0t}|^p + \io |\na u_0|^p + \io |\na u_0|^{p+2} + \io |\na \Theta_0|^p \le \eta^p
  \ee
  and
  for some $\eps\in (0,\eps_\star)$ and some
  $\eta\in (0,\eta_\star(\eps))$, then (\ref{0}) admits a global classical solution $(u,\Theta)$, fulfilling
  (\ref{tl1}) and (\ref{tl2}) with $\tm=\infty$ as well as $\Theta\ge 0$ in $\Om\times (0,\infty)$, which is such that moreover
  \be{12.6}
	\|\Theta(\cdot,t)-\Theta_\star\|_{L^\infty(\Om)} \le 2\eps
	\qquad \mbox{for all } t>0
  \ee
  and
  \be{12.5}
	\io |\na u_t(\cdot,t)|^p  + \io |\na u(\cdot,t)|^p + \io |\na\Theta(\cdot,t)|^p \le C\eta^p e^{-\kappa t}
	\qquad \mbox{for all } t>0
  \ee
  as well as
  \be{12.7}
	\int_0^\infty \io |\na (u_t+au)|^{p-2} |D^2(u_t+au)|^2 <\infty
  \ee
  and
  \be{12.8}
	\int_0^\infty \io |\na\Theta|^{p-2} |D^2\Theta|^2 <\infty.
  \ee
\end{theo}
\mysection{Local solvabilty. Proof of Theorem \ref{theo_loc}}\label{sect2}
In order to appropriately accentuate the dissipative features related to the third-order contribution to (\ref{0}),
formally substituting $v:=u_t+au$ we can rewrite (\ref{0}) in equivalent form as the parabolic-ODE-parabolic problem
\be{0v}
	\lball
	v_t = \na\cdot (\gamma(\Theta)\na v) + av - a^2 u + \na \cdot f(\Theta),
	\qquad & x\in\Om, \ t>0, \\[1mm]
	u_t = v-au,
	\qquad & x\in\Om, \ t>0, \\[1mm]
	\Theta_t = D\Del\Theta + \Gamma(\Theta) |\na v - a\na u|^2 + F(\Theta)\cdot (\na v - a\na u),
	\qquad & x\in\Om, \ t>0, \\[1mm]
	\frac{\pa v}{\pa\nu}=\frac{\pa u}{\pa\nu}=\frac{\pa\Theta}{\pa\nu}=0,
	\qquad & x\in\pO, \ t>0, \\[1mm]
	v(x,0)=v_0(x), \quad u(x,0)=u_0(x), \quad \Theta(x,0)=\Theta_0(x),
	\qquad & x\in\Om,
	\ear
\ee
which will be our actual object of study throughout the sequel.\abs
Spelt out this way, the problem becomes accessible to a fixed point based approach, which by suitably utilizing
scalar parabolic theory leads to the following statement on local existence and extensibility that can be viewed as
a rather verbatim translation of Theorem \ref{theo_loc} to the framework of (\ref{0v}):
\begin{lem}\label{lem_loc}
  Let $n\ge 1$ and $\Om\subset\R^n$ be a bounded domain with smooth boundary, and assume (\ref{gf}).
  Then whenever $\al\in (0,1)$ and
  \be{init}
	\lbal
	v_0\in C^{1+\al}(\bom)
	\mbox{ such that $\frac{\pa v_0}{\pa\nu}=0$ on $\pO$,} \\[1mm]
	u_0\in C^2(\bom)
	\mbox{ such that $\frac{\pa u_0}{\pa\nu}=0$ on $\pO$ \qquad and} \\[1mm]
	\Theta_0\in C^{1+\al}(\bom)
	\mbox{ such that $\Theta_0\ge 0$ in $\Om$ and $\frac{\pa \Theta_0}{\pa\nu}=0$ on $\pO$,}
	\ear
  \ee
  one can find $\tm\in (0,\infty]$ as well as
  \be{l1}
	\lbal
	v\in \Big( \bigcup_{\beta\in (0,1)} C^{1+\beta,\frac{1+\beta}{2}}(\bom\times [0,\tm))\Big) \cap C^{2,1}(\bom\times (0,\tm)),
		\\[1mm]
	u\in \Big( \bigcup_{\beta\in (0,1)} C^{1+\beta,\frac{1+\beta}{2}}(\bom\times [0,\tm))\Big) \cap C^{2,1}(\bom\times (0,\tm))
		\qquad \mbox{and} \\[1mm]
	\Theta\in \Big( \bigcup_{\beta\in (0,1)} C^{1+\beta,\frac{1+\beta}{2}}(\bom\times [0,\tm))\Big)
		\cap C^{2,1}(\bom\times (0,\tm))
	\ear
  \ee
  such that $\Theta\ge 0$ in $\Om\times (0,\tm)$, and that $(v,u,\Theta)$ solves (\ref{0v}) classically
  in $\Om\times (0,\tm)$. Moreover, $\tm$ can be chosen so as to be such that
  \bea{ext}
	& & \hs{-15mm}
	\mbox{if $\tm<\infty$, \quad then \quad} \nn\\
	& & \hs{-6mm}
	\limsup_{t\nearrow\tm} \Big\{ \|v(\cdot,t)-au(\cdot,t)\|_{W^{1,p}(\Om)} + \|\Theta(\cdot,t)\|_{L^\infty(\Om)} \Big\}
		=\infty
	\mbox{\quad for all $p\ge 2$ such that $p>n$.}
  \eea
\end{lem}
\proof
  In order to cleanly prepare our argument, given $(v_0,u_0,\Theta_0)$ fulfilling (\ref{init}) let us set
  \bea{l2}
	K &:=&
	\|v_0\|_{L^\infty(\Om)}
	+ \|\na v_0\|_{L^\infty(\Om)}
	+ \|v_0\|_{C^{1+\al}(\bom)}
	+ \|u_0\|_{L^\infty(\Om)}
	+ \|\na u_0\|_{L^\infty(\Om)}  \nn\\
	& & + \|\Theta_0\|_{L^\infty(\Om)}
	+ \|\Theta_0\|_{C^{1+\al}(\bom)}
  \eea
  as well as
  \be{l3}
	c_1\equiv c_1(K):=\max_{\xi\in [0,K+1]} \Gamma(\xi),
	\quad
	c_2\equiv c_2(K):=\min_{\xi\in [0,K+1]} \gamma(\xi)
	\quad \mbox{and} \quad
	c_3\equiv c_3(K):=\max_{\xi\in [0,K+1]} |F(\xi)|,
  \ee
  noting that $c_2>0$ by continuity and positivity of $\gamma$ in $[0,\infty)$.
  We furthermore abbreviate
  \be{l4}
	c_4\equiv c_4(K):=2c_1 \cdot (K+1)^2 + 2c_1 a^2\cdot (2K+1)^2
	+ c_3 \cdot (K+1) + c_3 a \cdot (2K+1),
  \ee
  and employ standard theory of gradient regularity in scalar parabolic equations (\cite{lieberman}) to find
  $\al_1=\al_1(K)\in (0,\al]$ and $c_5=c_5(K)>0$ with the property that whenever $T\in (0,1]$ and
  $\vp\in C^0(\bom\times [0,T)) \cap C^{2,1}(\bom\times (0,T))$ is such that $\frac{\pa\vp}{\pa\nu}=0$ on $\pO\times [0,T)$ and
  \be{l5}
	\|\vp(\cdot,0)\|_{C^{1+\al}(\bom)} \le K
  \ee
  as well as
  \be{l6}
	\|\vp_t - D\Del\vp\|_{L^\infty(\Om)} \le c_4
	\qquad \mbox{for all } t\in (0,T),
  \ee
  it follows that actually $\vp\in C^{1+\al_1,\frac{1+\al_1}{2}}(\bom\times [0,T])$ with
  \be{l7}
	\|\vp\|_{C^{1+\al_1,\frac{1+\al_1}{2}}(\bom\times [0,T])} \le c_5.
  \ee
  As $\gamma$, $\gamma'$ and $f'$ are locally Lipschitz continuous in $[0,\infty)$ according to our assumptions, this particularly
  implies the existence of $\al_2=\al_2(K)\in (0,\al_1]$, $c_6=c_6(K)>0$ and $c_7=c_7(K)>0$ such that for any such $\vp$ we have
  \be{l8}
	\|\gamma(\vp)\|_{C^{\al_2,\frac{\al_2}{2}}(\bom\times [0,T])}
	+ \|\gamma'(\vp)\na\vp\|_{C^{\al_2,\frac{\al_2}{2}}(\bom\times [0,T])}
	\le c_6
  \ee
  and
  \be{l87}
	\|\na \cdot f(\vp)\|_{L^\infty(\Om\times (0,T))}
	\le c_7.
  \ee
  Likewise, the same token from \cite{lieberman} provides $\al_3=\al_3(K)\in (0,\al_2]$ and $c_8=c_8(K)>0$
  such that if $T\in (0,1]$,
  then any triple $(A_1,A_2,\psi)$ of functions
  $A_1\in C^{\al_2,\frac{\al_2}{2}}(\bom\times [0,T])$,
  $A_2\in C^{\al_2,\frac{\al_2}{2}}(\bom\times [0,T];\R^n)$  and
  $\psi\in C^0(\bom\times [0,T)) \cap C^{2,1}(\bom\times (0,T))$ fulfilling
  \be{l9}
	A_1\ge c_2
	\qquad \mbox{on } \Om\times (0,T)
  \ee
  and
  \be{l10}
	\|A_1\|_{C^{\al_2,\frac{\al_2}{2}}(\bom\times [0,T])}
	+ \|A_2\|_{C^{\al_2,\frac{\al_2}{2}}(\bom\times [0,T])}
	\le c_6
  \ee
  as well as $\frac{\pa\psi}{\pa\nu}=0$ on $\pO\times [0,T)$ and
  \be{l11}
	\|\psi(\cdot,0)\|_{C^{1+\al}(\bom)} \le K
  \ee
  and
  \be{l12}
	\|\psi_t - A_1\Del\psi - A_2\cdot \na\psi\|_{L^\infty(\Om)} \le
	a\cdot (K+1) + a^2 \cdot (2K+1)
	+ c_7
	\qquad \mbox{for all } t\in (0,T)
  \ee
  in fact satisfies $\psi\in C^{1+\al_3,\frac{1+\al_3}{2}}(\bom\times [0,T])$ with
  \be{l13}
	\|\psi\|_{C^{1+\al_3,\frac{1+\al_3}{2}}(\bom\times [0,T])} \le c_8;
  \ee
  in particular, there exist $\al_4=\al_4(K)>0$ and $c_9=c_9(K)>0$ such that for each $\psi$ with these properties we have
  \be{l143}
	\|\psi(\cdot,t)\|_{L^\infty(\Om)}
	\le \|\psi(\cdot,0)\|_{L^\infty(\Om)} + c_9 t^{\al_4}
	\qquad \mbox{for all } t\in (0,T)
  \ee
  and
  \be{l14}
	\|\na\psi(\cdot,t)\|_{L^\infty(\Om)}
	\le \|\na\psi(\cdot,0)\|_{L^\infty(\Om)} + c_9 t^{\al_4}
	\qquad \mbox{for all } t\in (0,T).
  \ee
  We now let
  \be{l15}
	T\equiv T(K) := \min \Bigg\{ 1 \, , \, \frac{1}{2c_4(K)} \, , \, \frac{1}{c_9^{1/\al_4}(K)} \Bigg\},
  \ee
  and fixing any $\beta\in (0,\al_3)$, on the closed convex subset
  \be{l16}
	S:= \Big\{ v\in X \ \Big| \ \|v(\cdot,t)\|_{L^\infty(\Om)} \le K+1 \mbox{ and }
		\|\na v(\cdot,t)\|_{L^\infty(\Om)} \le K+1 \mbox{ for all } t\in [0,T] \Big\}
  \ee
  of the Banach space
  \bas
	X:=C^{1+\beta,\frac{1+\beta}{2}}(\bom\times [0,T])
  \eas
  we define a mapping $\Phi$ as follows: Given $\ovv\in S$, we write
  \be{l17}
	u(x,t):=e^{-at} u_0(x) + \int_0^t e^{-a(t-s)} \ovv(\cdot,s) ds,
	\qquad x\in\bom, \ t\in [0,T],
  \ee
  and observe that clearly, by (\ref{init}), also $u\in X$ with
  \bea{l186}
	\|u(\cdot,t)\|_{L^\infty(\Om)}
	&\le& e^{-at} \|u_0\|_{L^\infty(\Om)}
	+ \int_0^t e^{-a(t-s)} \|\ovv(\cdot,s)\|_{L^\infty(\Om)} ds \nn\\
	&\le& K + \int_0^t (K+1) ds \nn\\
	&\le& 2K+1
	\qquad \mbox{for all } t\in (0,T)
  \eea
  and, similarly,
  \bea{l18}
	\|\na u(\cdot,t)\|_{L^\infty(\Om)}
	&\le& 2K+1
	\qquad \mbox{for all } t\in (0,T)
  \eea
  according to the nonnegativity of $a$, the inequality $T\le 1$, and the definitions (\ref{l2}) and (\ref{l16}) of $T$ and $S$.\abs
  Now the H\"older continuity of $|\na \ovv - a\na u|^2$
  particularly ensures that in line with standard parabolic theory (\cite{LSU})
  there exist $\wh{T}\in (0,T]$ and a uniquely determined
  \be{l187}
	\Theta\in C^0(\bom\times [0,\wh{T})) \cap C^{2,1}(\bom\times (0,\wh{T}))
  \ee
  such that
  \be{l188}
	\|\Theta(\cdot,t)\|_{L^\infty(\Om)} < K+1
	\qquad \mbox{for all } t\in [0,\wh{T}),
  \ee
  that the problem
  \be{l19}
	\lball
	\Theta_t = D\Del\Theta + \Gamma(|\Theta|) |\na v-a\na u|^2 + F\big(|\Theta|\big) \cdot (\na v-a\na u),
	\qquad & x\in\Om, \ t\in (0,\wh{T}), \\[1mm]
	\frac{\pa\Theta}{\pa\nu}=0,
	\qquad & x\in\pO, \ t\in (0,\wh{T}), \\[1mm]
	\Theta(x,0)=\Theta_0(x),
	\qquad & x\in\Om,
	\ear
  \ee
  is solved in the classical sense, and that
  \be{l20}
	\mbox{if $\wh{T}<T$, \quad then \quad}
	\|\Theta(\cdot,\wh{T})\|_{L^\infty(\Om)} = K+1.
  \ee
  Here due to the assumption that $F(0)=0$ as well as the nonnegativity of $\Gamma$ and $\Theta_0$,
  a comparison principle firstly guarantees that
  \be{l21}
	\Theta\ge 0
	\qquad \mbox{in } \Om\times (0,\wh{T}),
  \ee
  and since (\ref{l3}) together with (\ref{l16}), (\ref{l18}) and (\ref{l4}) implies that
  \bea{l22}
	\Theta_t - D\Del\Theta
	&\le& |\Theta_t-D\Del\Theta| \nn\\
	&\le& 2c_1 |\na \ovv|^2
	+ 2c_1 a^2 |\na u|^2
	+ c_3 |\na \ovv|
	+ c_3 a |\na u| \nn\\
	&\le& 2c_1 \cdot (K+1)^2 + 2c_1 a^2 \cdot (2K+1)^2 + c_3 \cdot (K+1) + c_3 a\cdot (2K+1) \nn\\[1mm]
	&=& c_4
	\qquad \mbox{in } \Om\times (0,\wh{T}),
  \eea
  a second comparison argument based on the inequality $\|\Theta_0\|_{L^\infty(\Om)} \le K$ warrants that
  \bas
	\Theta(x,t) \le K+c_4 t
	\qquad \mbox{for all } (x,t)\in\Om\times (0,\wh{T}).
  \eas
  As $c_4\wh{T} \le c_4 T \le \frac{1}{2}$ by (\ref{l15}), together with (\ref{l21}) and (\ref{l20}) this ensures that
  $\wh{T}$ cannot be smaller than $T$, and that hence the inequalities in
  \be{l23}
	0 \le \Theta \le K+1
  \ee
  and in (\ref{l22}) actually hold throughout $\Om\times (0,T)$.
  This enables us to rely on our considerations near (\ref{l7})-(\ref{l8}) to see that thanks to (\ref{l2}),
  \be{l221}
	\Theta\in C^{1+\al_1,\frac{1+\al_1}{2}}(\bom\times [0,T]),
  \ee
  with
  \bas
	\|\gamma(\Theta)\|_{C^{\al_2,\frac{\al_2}{2}}(\bom\times [0,T])}
	+ \|\gamma'(\Theta)\na\vp\|_{C^{\al_2,\frac{\al_2}{2}}(\bom\times [0,T])}
	\le c_5.
  \eas
  Since $A_1:=\gamma(\Theta)$ and $A_2:=\gamma'(\Theta) \na\Theta$ thus comply with (\ref{l10}), 
  combining (\ref{l3}) with
  (\ref{l23}) in verifying (\ref{l9}), recalling (\ref{l2}) in asserting (\ref{l11}), and deducing
  (\ref{l12}) on the basis of the estimate
  \bas
	\big\|a\ovv - a^2 u + \na\cdot f(\Theta)\big\|_{L^\infty(\Om)}
	\le a\|\ovv\|_{L^\infty(\Om)}
	+ a^2 \|u\|_{L^\infty(\Om)}
	+ c_7
	\le a\cdot (K+1) + a^2\cdot (2K+1) + c_7
  \eas
  which for all $t\in (0,T)$ is implied by (\ref{l16}), (\ref{l186}) and (\ref{l87}),
  we may now draw on
  (\ref{l13}), (\ref{l143}) and (\ref{l14}) to see that the unique solution
  \be{l222}
	v\in C^0(\bom\times [0,T)) \cap C^{2,1}(\bom\times (0,T))
  \ee
  of the linear and uniformly parabolic problem
  \be{l24}
	\lball
	v_t = \na\cdot (\gamma(\Theta)\na v) + a\ovv - a^2 u + \na \cdot f(\Theta),
	\qquad & x\in\Om, \ t\in (0,T), \\[1mm]
	\frac{\pa v}{\pa\nu}=0,
	\qquad & x\in\pO, \ t\in (0,T), \\[1mm]
	v(x,0)=v_0(x),
	\qquad & x\in\Om,
	\ear
  \ee
  the existence of which again being guaranteed by standard parabolic theory (\cite{LSU}), indeed satisfies
  \be{l25}
	v\in C^{1+\al_3,\frac{1+\al_3}{2}}(\bom\times [0,T])
	\qquad \mbox{with} \qquad
	\|v\|_{C^{1+\al_3,\frac{1+\al_3}{2}}(\bom\times [0,T])} \le c_8
  \ee
  as well as
  \be{l265}
	\|v(\cdot,t)\|_{L^\infty(\Om)}
	\le \|v_0\|_{L^\infty(\Om)} + c_9 t^{\al_4}
	\le K + c_9 T^{\al_4}
	\le K+1
	\qquad \mbox{for all } t\in [0,T]
  \ee
  and, similarly,
  \be{l26}
	\|\na v(\cdot,t)\|_{L^\infty(\Om)}
	\le K+1
	\qquad \mbox{for all } t\in [0,T]
  \ee
  thanks to (\ref{l2}) and (\ref{l15}).\abs
  Consequently, if now we let $\Phi \ovv := v$, then since $\beta\le \al_3$ it follows from (\ref{l25}) that
  $\Phi \ovv \in X$, whereupon (\ref{l26}) implies that, in fact, $\Phi\ovv \in S$.
  As the strictness of the inequality $\beta<\al_3$ ensures compactness of the embedding
  $C^{1+\al_3,\frac{1+\al_3}{2}}(\bom\times [0,T]) \hra X$, from (\ref{l25}) we thus infer that $\Phi (S)$
  is a compact subset of $S$.
  Since furthermore a straightforward argument based on the uniqueness properties of (\ref{l19}) and (\ref{l24}) shows
  that if $(\ovv_j)_{j\in\N} \subset S$ and $\ovv\in S$ are such that $\ovv_j\to \ovv$ in $X$, then
  $(\Phi \ovv_j)_{j\in\N}$ cannot have any accumulation point in $X$ other than $\Phi \ovv$,
  we finally obtain that $\Phi$ is also continuous on $S$, whence the Schauder fixed point theorem ensures the existence
  of an element $v$ of $S$ fulfilling $\Phi v=v$.\abs
  In line with our construction of $S$ and the inclusions in (\ref{l222}) and (\ref{l221}),
  along with the functions $u$ and $\Theta$ accordingly defined through (\ref{l17}) and (\ref{l19}) this fixed point forms a
  classical solution of (\ref{0v}) in $\Om\times (0,T)$, with regularity properties consistent with those listed in (\ref{l1}),
  and with nonnegativity of $\Theta$ implied by (\ref{l23}).\abs
  Finally, according to a standard prolongation argument, the obtained solution can be extended up to some maximal
  $\tm\in (0,\infty]$ which satisfies
  \be{e1}
	\mbox{if $\tm<\infty$, \quad then \quad}
	\limsup_{t\nearrow\tm} \Big\{
	 \|v(\cdot,t)\|_{C^{1+\beta}(\bom)}
	 +\|u(\cdot,t)\|_{W^ {1,\infty}(\Om)}
	 +\|\Theta(\cdot,t)\|_{C^{1+\beta}(\bom)}
	\Big\}
	=\infty,
  \ee
  and in order to make sure that in fact even (\ref{ext}) holds, let us suppose that $\tm<\infty$, but that with some
  $c_{10}>0$ and some
  $p\ge 2$ such that $p>n$ we had
  \be{e2}
	\|v(\cdot,t)-a u(\cdot,t)\|_{W^{1,p}(\Om)} + \|\Theta(\cdot,t)\|_{L^\infty(\Om)} \le c_{10}
	\qquad \mbox{for all } t\in (0,\tm).
  \ee
  Then, in particular, $\Gamma(\Theta)\in L^\infty(\Om\times (0,\tm))$, so that from (\ref{0v}) and (\ref{e2}) we would obtain
  $c_{11}>0$ such that
  \bas
	\|\Theta_t - D\Del\Theta\|_{L^\frac{p}{2}(\Om)} \le c_{11}
	\qquad \mbox{for all } t\in (0,\tm).
  \eas
  According to the inequalities $\frac{p}{2}\ge 1$ and $\frac{p}{2}>\frac{n}{2}$, due to a known H\"older regularity result
  (\cite{PV}), and again due to the local Lipschitz continuity of $\gamma$ and $f$,
  this would entail the existence of $\al_5\in (0,1)$ such that
  \be{e3}
	\Theta\in C^{\al_5,\frac{\al_5}{2}}(\bom\times [0,\tm])
  \ee
  as well as
  \bas
	\big\{ \gamma(\Theta) \, , \, f(\Theta) \big\} \subset C^{\al_5,\frac{\al_5}{2}}(\bom\times [0,\tm]),
  \eas
  whereupon \cite{lieberman} would apply to the first equation in (\ref{0v}), here interpreted as the scalar identity
  \bas
	v_t=\na \cdot \big(\gamma(\Theta)\na v + f(\Theta)\big) + h(x,t),
	\qquad (x,t)\in \Om\times (0,\tm),
  \eas
  with $h:= av-a^2 u$,
  so as to yield $\al_6\in (0,1)$ such that
  \be{e33}
	v\in C^{1+\al_6,\frac{1+\al_6}{2}}(\bom\times [0,\tm]),
  \ee
  because for some $c_{12}>0$ and $c_{13}>0$,
  \be{e32}
	c_{12} \le \gamma(\Theta)\le c_{13}
	\qquad \mbox{in } \Om\times (0,\tm)
  \ee
  thanks to (\ref{e2}) and the positivity and local boundedness of $\gamma$ on $[0,\infty)$,
  and because (\ref{e2}) would moreover entail that $h\in L^\infty(\Om\times (0,\tm))$.
  In view of the second equation in (\ref{0v}), this would mean that also
  \be{e34}
	u\in C^{1+\al_6,\frac{1+\al_6}{2}}(\bom\times [0,\tm]),
  \ee
  whence actually $\na v-a\na u \in C^{\al_7,\frac{\al_7}{2}}(\bom\times [0,\tm])$ for some $\al_7\in (0,1)$.
  Together with the fact that for suitably small $\al_8\in (0,1)$ we have
  \be{e4}
	\Gamma(\Theta) \in C^{\al_8,\frac{\al_8}{2}}(\bom\times [0,\tm])
  \ee
  by (\ref{e3}) and the local Lipschitz continuity of $\Gamma$ in $[0,\infty)$, this would enable us to apply classical interior
  parabolic Schauder theory (\cite{LSU}) to the third equation in (\ref{0v}) to infer the existence of $\al_9\in (0,1)$ such that
  \be{e5}
	\Theta\in C^{2+\al_9,1+\frac{\al_9}{2}} \mbox{$(\bom\times [\frac{1}{4}\tm,\tm])$}.
  \ee
  Thus,
  \bas
	v_t=\gamma(\Theta)\Del v + \wh{h}(x,t),
	\qquad (x,t)\in \Om\times (0,\tm),
  \eas
  with $\wh{h}:=\gamma'(\Theta)\na\theta\cdot\na v + av-a^2 u+\na\cdot f(\Theta)$ satisfying
  $\wh{f}\in C^{\al_{10},\frac{\al_{10}}{2}}(\bom\times [\frac{1}{4}\tm,\tm])$ for some $\al_{10}\in (0,1)$
  due to (\ref{e5}), (\ref{e33}), (\ref{e34}) and the Lipschitz continuity of $\gamma'$ and $f'$
  on $[0,\|\Theta\|_{L^\infty(\Om\times (0,\tm))}]$, so that again by (\ref{e4}) and (\ref{e32}), another application of
  parabolic Schauder theory would yield $\al_{11}\in (0,1)$ such that
  \be{e6}
	v\in C^{2+\al_{11},1+\frac{\al_{11}}{2}}(\bom\times [\mbox{$\frac{1}{2}$}\tm,\tm]).
  \ee
  Once more using that $u_t=v-au$, from (\ref{e6}) and the inclusion $u(\cdot,\frac{1}{2}\tm) \in C^2(\bom)$ we would obtain
  $c_{13}>0$ such that
  \bas
	\|u(\cdot,t)\|_{C^2(\bom)} \le c_{13}
	\qquad \mbox{for all } t\in \Big(\frac{1}{2}\tm,\tm\Big),
  \eas
  which in conjunction with (\ref{e6}) and (\ref{e5}) would contradict (\ref{e1}).
\qed
Our result on local-in-time solvability and extensibility in (\ref{0}) has thereby in fact been established already:\abs
\proofc of Theorem \ref{theo_loc}.\quad
  We only need to let $\tm$ and $(v,u,\Theta)$ be as provided by Lemma \ref{lem_loc} when applied to
  $(v_0,u_0,\Theta_0):=(u_{0t}+au_0,u_0,\Theta_0)$,
  and to use (\ref{l1}), (\ref{0v}) and (\ref{ext})
  to verify (\ref{tl1}), (\ref{tl2}), (\ref{0}) and (\ref{Ext}).
\qed
\mysection{Small-data solutions. Proof of Theorem \ref{theo12}}\label{sect3}
\subsection{Basic evolution features of gradients}
Our approach toward the derivation of the result on global existence of small-data solutions formulated in
Theorem \ref{theo12} will be launched by three basic observations on elementary evolution features enjoyed by
gradients of the solution components in (\ref{0v}).
The first of these concerns the quantity $v$.
\begin{lem}\label{lem3}
  Let $n\ge 1$ and $p\ge 2$, and assume (\ref{gf}) and (\ref{init}). Then
  \bea{3.1}
	\hs{-8mm}
	\frac{1}{p} \frac{d}{dt} \io |\na v|^p
	+ \io \gamma(\Theta) |\na v|^{p-2} |D^2 v|^2
	&\le& - \frac{1}{2} \io \gamma'(\Theta) |\na v|^{p-2} \na\Theta\cdot\na |\na v|^2 \nn\\
	& & - \io \gamma'(\Theta) (\na\Theta\cdot\na v) \na v\cdot\na |\na v|^{p-2} \nn\\
	& & + 2a\io |\na v|^p
	+ a^{p+1} \io |\na u|^p \nn\\
	& & - \io \big( f'(\Theta)\cdot\na\Theta\big) |\na v|^{p-2} \Del v \nn\\
	& & - \io \big( f'(\Theta)\cdot\na\Theta\big) \na v \cdot \na |\na v|^{p-2} \nn\\
	& & + \frac{1}{2} \int_{\pO} \gamma(\Theta) |\na v|^{p-2} \frac{\pa |\na v|^2}{\pa\nu}
	\qquad \mbox{for all } t\in (0,\tm).
  \eea
\end{lem}
\proof
  We use the first equation in (\ref{0v}) and integrate by parts to compute
  \bea{3.2}
	\frac{1}{p} \frac{d}{dt} \io |\na v|^p
	&=& \io |\na v|^{p-2} \na v\cdot\na \Big\{ \gamma(\Theta) \Del v + \gamma'(\Theta) \na\Theta\cdot\na v + av-a^2 u
		+ \na\cdot f(\Theta) \Big\}
		\nn\\
	&=& \io |\na v|^{p-2} \na v\cdot \Big\{ \gamma(\Theta) \na\Del v + \gamma'(\Theta) \Del v \na\Theta \Big\} \nn\\
	& & - \io \Big\{ |\na v|^{p-2} \Del v + \na |\na v|^{p-2} \cdot\na v \Big\} \gamma'(\Theta) \na\Theta\cdot\na v \nn\\
	& & + a\io |\na v|^p - a^2 \io |\na v|^{p-2} \na v\cdot\na u \nn\\
	& & - \io \Big\{ |\na v|^{p-2} \Del v + \na |\na v|^{p-2} \cdot\na v \Big\} f'(\Theta)\cdot\na\Theta
	\qquad \mbox{for all } t\in (0,\tm),
  \eea
  where since $\na v\cdot\na \Del v=\frac{1}{2} \Del |\na v|^2 - |D^2 v|^2$,
  \bas
	\io |\na v|^{p-2} \na v\cdot \gamma(\Theta) \na\Del v
	&=& \frac{1}{2} \io \gamma(\Theta) |\na v|^{p-2} \Del |\na v|^2
	- \io \gamma(\Theta) |\na v|^{p-2} |D^2 v|^2 \\
	&=& - \frac{4(p-2)}{p^2} \io \gamma(\Theta) \Big| \na |\na v|^\frac{p}{2}\Big|^2
	- \frac{1}{2} \io \gamma'(\Theta) |\na v|^{p-2} \na\Theta\cdot\na |\na v|^2 \\
	& & + \frac{1}{2} \int_{\pO} \gamma(\Theta) |\na v|^{p-2} \frac{\pa |\na v|^2}{\pa\nu}
	- \io \gamma(\Theta) |\na v|^{p-2} |D^2 v|^2
  \eas
  for all $t\in (0,\tm)$.
  Since Young's inequality warrants that
  \bas
	- a^2 \io |\na v|^{p-2} \na v\cdot\na u
	&\le& a^2 \io |\na v|^{p-1} |\na u| \\
	&=& \io \big\{ a|\na v|^p \big\}^\frac{p-1}{p} \cdot a^\frac{p+1}{p} |\na u| \\
	&\le& a \io |\na v|^p
	+ a^{p+1} \io |\na u|^p
	\qquad \mbox{for all } t\in (0,\tm),
  \eas
  due to a cancellation of two summands in (\ref{3.2}) containing $\Del v$ we thus obtain that
  for all $t\in (0,\tm)$,
  \bas
	\frac{1}{p} \frac{d}{dt} \io |\na v|^p
	&\le& - \io \gamma(\Theta) |\na v|^{p-2} |D^2 v|^2
	- \frac{4(p-2)}{p^2} \io \gamma(\Theta) \Big| \na |\na v|^\frac{p}{2}\Big|^2 \\
	& & - \frac{1}{2} \io \gamma'(\Theta) |\na v|^{p-2} \na\Theta\cdot\na |\na v|^2
	- \io \gamma'(\Theta) (\na\Theta\cdot\na v) \na v\cdot\na |\na v|^{p-2} \\
	& & + 2a\io |\na v|^p
	+ a^{p+1} \io |\na u|^p \\
	& & - \io \big( f'(\Theta)\cdot\na\Theta\big) |\na v|^{p-2} \Del v
	- \io \big( f'(\Theta)\cdot\na\Theta\big) \na v \cdot \na |\na v|^{p-2} \\
	& & + \frac{1}{2} \int_{\pO} \gamma(\Theta) |\na v|^{p-2} \frac{\pa |\na v|^2}{\pa\nu},
  \eas
  from which (\ref{3.1}) follows in view of our assumption that $p\ge 2$.
\qed
Due to the simple ODE structure determining their evolution, gradients of the second solution component can be described
in a comparatively simple manner:
\begin{lem}\label{lem4}
  If $n\ge 1$ and $p\ge 2$, and if (\ref{gf}) and (\ref{init}) hold, then
  \be{4.1}
	\frac{1}{p} \frac{d}{dt} \io |\na u|^p
	+ \frac{a}{2} \io |\na u|^p
	\le \Big(\frac{2}{a}\Big)^{p-1} \io |\na v|^p
	\qquad \mbox{for all } t\in (0,\tm)
  \ee
  and
  \be{4.2}
	\frac{1}{p+2} \frac{d}{dt} \io |\na u|^{p+2}
	+ \frac{a}{2} \io |\na u|^{p+2}
	\le \Big(\frac{2}{a}\Big)^{p+1} \io |\na v|^{p+2}
	\qquad \mbox{for all } t\in (0,\tm).
  \ee
\end{lem}
\proof
  By means of Young's inequality, using the second equation in (\ref{0v}) we can estimate
  \bas
	\frac{1}{p} \frac{d}{dt} \io |\na u|^p
	&=& \io |\na u|^{p-2} \na v\cdot\na u
	- a \io |\na u|^p \\
	&\le& \io |\na u|^{p-1} |\na v|
	- a \io |\na u|^p \\
	&=& \io \Big\{ \frac{a}{2} |\na u|^p \Big\}^\frac{p-1}{p} \cdot \Big(\frac{2}{a}\Big)^\frac{p-1}{p} |\na v|
	- a \io |\na u|^p \\
	&\le& \frac{a}{2} \io |\na u|^p
	+ \Big(\frac{2}{a}\Big)^{p-1} \io |\na v|^p
	- a \io |\na u|^p
	\qquad \mbox{for all } t\in (0,\tm),
  \eas
  and thereby obtain (\ref{4.1}). Replacing $p$ with $p+2$ here directly yields (\ref{4.2}).
\qed
The temperature distribution, finally, will again be influenced by diffusion, as is reflected in the following inequality
describing its gradient:
\begin{lem}\label{lem5}
  Let $n\ge 1$. Then
  for all $p\ge 2$ there exists $K_0=K_0(p)>0$ such that if (\ref{gf}) and (\ref{init}) hold, we have
  \bea{5.1}
	& & \hs{-30mm}
	\frac{1}{p} \frac{d}{dt} \io |\na \Theta|^p
	+ \frac{D}{2} \io |\na\Theta|^{p-2} |D^2\Theta|^2 \nn\\
	&\le& \frac{K_0}{D} \io \Gamma^2(\Theta) |\na \Theta|^{p-2} |\na v|^4
	+ \frac{K_0 a^4}{D} \io \Gamma^2(\Theta) |\na \Theta|^{p-2} |\na u|^4 \nn\\
	& & + \frac{K_0}{D} \io \big|F(\Theta)|^2 |\na\Theta|^{p-2} |\na v|^2
	+ \frac{K_0 a^2}{D} \io \big|F(\Theta)|^2 |\na\Theta|^{p-2} |\na u|^2 \nn\\
	& & + \frac{D}{2} \int_{\pO} |\na\Theta|^{p-2} \frac{\pa |\na\Theta|^2}{\pa\nu}
	\qquad \mbox{for all } t\in (0,\tm).
  \eea
\end{lem}
\proof
  Thanks to the identity $\na\Theta\cdot\na\Del\Theta=\frac{1}{2}\Del |\na \Theta|^2 - |D^2\Theta|^2$, an integration by parts
  in the third equation from (\ref{0v}) shows that
  \bas
	\frac{1}{p} \frac{d}{dt} \io |\na\Theta|^p
	&=& \io |\na\Theta|^{p-2} \na\Theta\cdot \na \Big\{ D\Del\Theta + \Gamma(\Theta) |\na v-a\na u|^2
	+ F(\Theta)\cdot (\na v-a\na u) \Big\} \\
	&=& \frac{D}{2} \io |\na\Theta|^{p-2} \Del |\na\Theta|^2
	- D \io |\na\Theta|^{p-2} |D^2\Theta|^2 \\
	& & - \io \Gamma(\Theta) \cdot \Big\{ \na\Theta\cdot\na |\na\Theta|^{p-2} + |\na\Theta|^{p-2} \Del\Theta\Big\} \cdot
		|\na v-a\na u|^2 \\
	& & - \io \big\{ F(\Theta) \cdot (\na v-a\na u)\big\}
		\Big\{ \na\Theta\cdot\na |\na\Theta|^{p-2} + |\na\Theta|^{p-2} \Del\Theta\Big\} \\
	&\le& \frac{D}{2} \int_{\pO} |\na\Theta|^{p-2} \frac{\pa |\na\Theta|^2}{\pa\nu}
	- D \io |\na\Theta|^{p-2} |D^2\Theta|^2 \\
	& & + (p-2+\sqrt{n}) \io \Gamma(\Theta) |\na\Theta|^{p-2} |D^2\Theta| \cdot |\na v-a\na u|^2 \\
	& & + (p-2+\sqrt{n}) \io \big| F(\Theta)\big| |\na\Theta|^{p-2} |D^2\Theta| |\na v-a\na u|
	\qquad \mbox{for all } t\in (0,\tm),
  \eas
  because $\na |\na\Theta|^{p-2} = 
  (p-2) |\na\Theta|^{p-4} D^2\Theta\cdot\na\Theta$ in $\{\na\Theta\ne 0\}$,
  because thus $\na |\na\Theta|^{p-2} \cdot \na |\na\Theta|^2 \ge 0$, and because $|\Del\Theta| \le \sqrt{n} |D^2\Theta|$.
  As Young's inequality implies that
  \bas
	& & \hs{-15mm}
	(p-2+\sqrt{n}) \io \Gamma(\Theta) |\na\Theta|^{p-2} |D^2\Theta| \cdot |\na v-a\na u|^2 \\
	&\le& \frac{D}{4} \io |\na\Theta|^{p-2} |D^2\Theta|^2
	+ \frac{(p-2+\sqrt{n})^2}{D} \io \Gamma^2(\Theta) |\na\Theta|^{p-2} |\na v-a\na u|^4 \\
	&\le& \frac{D}{4} \io |\na\Theta|^{p-2} |D^2\Theta|^2 \\
	& & + \frac{8(p-2+\sqrt{n})^2}{D} \io \Gamma^2(\Theta) |\na\Theta|^{p-2} |\na v|^4
	+ \frac{8(p-2+\sqrt{n})^2 a^4}{D} \io \Gamma^2(\Theta) |\na\Theta|^{p-2} |\na u|^4
  \eas
  and
  \bas
	& & \hs{-15mm}
	(p-2+\sqrt{n}) \io \big| F(\Theta)\big| |\na\Theta|^{p-2} |D^2\Theta| |\na v-a\na u| \\
	&\le& \frac{D}{4} \io |\na\Theta|^{p-2} |D^2\Theta|^2
	+ \frac{(p-2+\sqrt{n})^2}{D} \io \big| F(\Theta)\big|^2 |\na\Theta|^{p-2} |\na v-a\na u|^2 \\
	&\le& \frac{D}{4} \io |\na\Theta|^{p-2} |D^2\Theta|^2 \\
	& & + \frac{2(p-2+\sqrt{n})^2}{D} \io \big| F(\Theta)\big|^2 |\na\Theta|^{p-2} |\na v|^2
	+ \frac{2(p-2+\sqrt{n})^2 a^2}{D} \io \big| F(\Theta)\big|^2 |\na\Theta|^{p-2} |\na u|^2
  \eas
  for all $t\in (0,\tm)$,
  this entails (\ref{5.1}) if we choose $K_0(p)>0$ suitably large.
\qed
\subsection{Linear dissipation dominates superlinear forces along small-data trajectories}
In order to prepare an appropriate exploitation of the second-order dissipative contributions to (\ref{3.1}) and (\ref{5.1}),
let us state and derive a consequence of a standard Poincar\'e type inequality.
\begin{lem}\label{lem6}
  Let $n\ge 1$, and for fixed $p\ge 2$ let $C_P(p)>0$ be such that
  \be{6.1}
	\io \bigg| \vp - \frac{1}{|\Om|} \io \vp \bigg|^p \le C_P(p) \io |\na\vp|^p
	\qquad \mbox{for all } \vp\in W^{1,p}(\Om).
  \ee
  Then
  \be{6.2}
	\io |\na\vp|^p \le (p-2+\sqrt{n})^2 C_P^\frac{2}{p}(p) \io |\na\vp|^{p-2} |D^2\vp|^2
	\qquad \mbox{for all $\vp\in C^2(\bom)$ fulfilling $\frac{\pa\vp}{\pa\nu}=0$ on } \pO.
  \ee
\end{lem}
\proof
  For $\eta>0$, writing $\ov{\vp}:=\frac{1}{|\Om|} \io \vp$ and using that $\frac{\pa\vp}{\pa\nu}=0$ on $\pO$
  we integrate by parts in confirming that
  \bea{6.3}
	\io \big(|\na\vp|^2+\eta\big)^\frac{p}{2}
	&=& \io \big(|\na\vp|^2+\eta\big)^\frac{p-2}{2} \big(|\na\vp|^2+\eta\big) \nn\\
	&=& \io \big(|\na\vp|^2+\eta\big)^\frac{p-2}{2} \na\vp\cdot\na (\vp-\ov{\vp})
	+ \eta \io \big(|\na\vp|^2+\eta\big)^\frac{p-2}{2} \nn\\
	&=& - \io (\vp-\ov{\vp}) \na\cdot \Big\{ \big(|\na\vp|^2+\eta\big)^\frac{p-2}{2} \na\vp\Big\}
	+ \eta \io \big(|\na\vp|^2+\eta\big)^\frac{p-2}{2}.
  \eea
  Since
  \bas
	\bigg| \na \cdot \Big\{ \big(|\na\vp|^2+\eta\big)^\frac{p-2}{2} \na\vp\Big\} \bigg|
	&=& \Big| (p-2) \big(|\na\vp|^2+\eta\big)^\frac{p-4}{2} \na\vp\cdot (D^2\vp\cdot\na\vp)
	+ \big(|\na\vp|^2+\eta\big)^\frac{p-2}{2} \Del\vp \Big| \\
	&\le& (p-2) \big(|\na\vp|^2+\eta\big)^\frac{p-4}{2} |\na\vp|^2 |D^2\vp|
	+ \sqrt{n} \big(|\na\vp|^2+\eta\big)^\frac{p-2}{2} |D^2\vp| \\
	&\le& (p-2+\sqrt{n}) \big(|\na\vp|^2+\eta\big)^\frac{p-2}{2} |D^2\vp|
	\qquad \mbox{in } \Om
  \eas
  due to the inequality $|\Del\vp|\le\sqrt{n} |D^2\vp|$, we can here combine the H\"older inequality with (\ref{6.1}) to see that
  \bas
	& & \hs{-14mm}
	- \io (\vp-\ov{\vp}) \na\cdot \Big\{ \big(|\na\vp|^2+\eta\big)^\frac{p-2}{2} \na\vp\Big\} \\
	&\le& (p-2+\sqrt{n}) \io |\vp-\ov{\vp}| \big(|\na\vp|^2+\eta\big)^\frac{p-2}{2} |D^2\vp| \\
	&\le& (p-2+\sqrt{n}) \cdot \bigg\{ \io (\vp-\ov{\vp})^2 \big(|\na\vp|^2+\eta\big)^\frac{p-2}{2} \bigg\}^\frac{1}{2} \cdot
		\bigg\{ \io \big(|\na\vp|^2+\eta\big)^\frac{p-2}{2} |D^2\vp|^2 \bigg\}^\frac{1}{2} \\
	&\le& (p-2+\sqrt{n}) \cdot \bigg\{ \io |\vp-\ov{\vp}|^p \bigg\}^\frac{1}{p} \cdot
		\bigg\{ \io \big(|\na\vp|^2+\eta\big)^\frac{p}{2} \bigg\}^\frac{p-2}{2p} \cdot
		\bigg\{ \io \big(|\na\vp|^2+\eta\big)^\frac{p-2}{2} |D^2\vp|^2 \bigg\}^\frac{1}{2} \\
	&\le& (p-2+\sqrt{n}) C_P^\frac{1}{p}(p) \cdot \bigg\{ \io |\na\vp|^p \bigg\}^\frac{1}{p} \cdot
		\bigg\{ \io \big(|\na\vp|^2+\eta\big)^\frac{p}{2} \bigg\}^\frac{p-2}{2p} \cdot
		\bigg\{ \io \big(|\na\vp|^2+\eta\big)^\frac{p-2}{2} |D^2\vp|^2 \bigg\}^\frac{1}{2}.
  \eas
  As $\io |\na\vp|^p \le \io \big(|\na\vp|^2+\eta\big)^\frac{p}{2}$, we may thus invoke Young's inequality to estimate
  \bas
	& & \hs{-24mm}
	- \io (\vp-\ov{\vp}) \na\cdot \Big\{ \big(|\na\vp|^2+\eta\big)^\frac{p-2}{2} \na\vp\Big\} \\
	&\le& (p-2+\sqrt{n}) C_P^\frac{1}{p}(p) \cdot \bigg\{ \io \big(|\na\vp|^2+\eta\big)^\frac{p}{2} \bigg\}^\frac{1}{2} \cdot
		\bigg\{ \io \big(|\na\vp|^2+\eta\big)^\frac{p-2}{2} |D^2\vp|^2 \bigg\}^\frac{1}{2} \\
	&\le& \frac{1}{2} \io \big(|\na\vp|^2+\eta\big)^\frac{p}{2}
	+ \frac{1}{2} (p-2+\sqrt{n})^2 C_P^\frac{2}{p}(p) \io \big(|\na\vp|^2+\eta\big)^\frac{p-2}{2} |D^2\vp|^2,
  \eas
  so that (\ref{6.3}) shows that for all $\eta>0$,
  \bas
	\io \big(|\na\vp|^2+\eta\big)^\frac{p}{2}
	\le (p-2+\sqrt{n})^2 C_P^\frac{2}{p}(p) \io \big(|\na\vp|^2+\eta\big)^\frac{p-2}{2} |D^2\vp|^2
	+ 2\eta \io \big(|\na\vp|^2+\eta\big)^\frac{p-2}{2}.
  \eas
  Taking $\eta\searrow 0$ yields (\ref{6.2}).
\qed
By making use of this functional inequality,
we can turn (\ref{3.1}) and (\ref{4.1}) into a superlinearly forced but linearly damped ODI for
some linear combination of $t\mapsto \io |\na v|^p$ and $t\mapsto \io |\na u|^p$,
provided that the parameter $a$ is small relative to $\gamma(\Theta)$, and
provided that moreover the rightmost summand in (\ref{3.1}) is favorably signed.
As a convenient sufficient assumption guaranteeing the latter, in the following lemma and its descendants we will
require $\Om$ to be convex.
\begin{lem}\label{lem7}
  Let $n\ge 1$, and suppose that $\Om$ is convex.
  Then for all $p\ge 2$ one can choose $\del_1(p)>0, k_1=k_1(p)>0$ and $K_1=K_1(p)>0$ in such a way that if
  beyond (\ref{gf}) we have
  \be{7.1}
	\frac{a}{\gamma(\Theta_\star)} \le \del_1(p)
  \ee
  with some $\Theta_\star\ge 0$,
  one can find $\eps_1=\eps_1(p,a,D,\gamma,\Gamma,f,F,\Theta_\star)>0$
  such that whenever (\ref{init}) holds and $T\in (0,\tm]$ is such that
  \be{7.2}
	|\Theta-\Theta_\star| \le \eps
	\qquad \mbox{in } \Om\times (0,T)
  \ee
  with some $\eps\in (0,\eps_1)$, it follows that
  for all $t\in (0,T)$,
  \bea{7.3}
	& & \hs{-16mm}
	\frac{d}{dt} \bigg\{ \io |\na v|^p + 8a^{p-1} \gamma_\star \del_1(p) \io |\na u|^p \bigg\}
	+ k_1 \gamma_\star \io |\na v|^{p-2} |D^2 v|^2
	+ k_1 \gamma_\star \io |\na v|^p
	+ k_1 a^p \gamma_\star \io |\na u|^p \nn\\
	&\le& K_1 \gamma_\star^{1-p} \|f'\|_{L^\infty(I_1)}^p \io |\na\Theta|^p
	+ K_1 \gamma_\star^{-\frac{p+2}{2}} \|\gamma'\|_{L^\infty(I_1)}^{p+2} \io |\na\Theta|^{p+2}
	+ K_1 \io |\na v|^{p+2}
  \eea
  where
  \be{7.44}
	I_1:=\big\{ \xi\ge 0 \ | \ |\xi-\Theta_\star| \le \eps_1 \big\}
  \ee
  and
  \be{7.4}
	\gamma_\star:=\inf_{\xi\in I_1} \gamma(\xi).
  \ee
\end{lem}
\proof
  According to Lemma \ref{lem6}, there exists $c_1=c_1(p)>0$ such that
  \be{7.5}
	\io |\na\vp|^p
	\le c_1 \io |\na\vp|^{p-2} |D^2\vp|^2
	\qquad \mbox{for all $\vp\in C^2(\bom)$ satisfying $\frac{\pa\vp}{\pa\nu}=0$ on $\pO$,}
  \ee
  and we let
  \be{7.6}
	\del_1\equiv\del_1(p)
	:= \frac{1}{32(1+2^p)c_1}.
  \ee
  Henceforth fixing any $a>0$, $D>0$ and $\Theta_\star\ge 0$
  as well as $0<\gamma\in C^2([0,\infty))$, $0\le\Gamma\in C^1([0,\infty))$,
  $f\in C^2([0,\infty);\R^n)$ and $F\in C^1([0,\infty);\R^n)$ such that $F(0)=0$
  and that (\ref{7.1}) holds, by continuity of
  $\gamma$ we can pick $\eps_1=\eps_1(p,a,D,\gamma,\Gamma,F,f,\Theta_\star)>0$ such that
  taking $I_1$ and $\gamma_\star$ as in (\ref{7.44}) and (\ref{7.4}) we have
  \bas
	\gamma_\star \ge \frac{\gamma(\Theta_\star)}{2}
  \eas
  and thus, in particular,
  \be{7.7}
	a \le 2\del_1 \gamma_\star,
  \ee
  and we thereupon assume that (\ref{init}) is satisfied, and that $T\in (0,\tm]$ is
  such that (\ref{7.2}) is valid with some $\eps\in (0,\eps_1)$.\abs
  Then using that $\frac{\pa |\na v|^2}{\pa\nu} \le 0$ on $\pO\times (0,T)$ by convexity of $\Om$ (\cite{lions_ARMA}),
  abbreviating $\gamma_{\star\star}:=\|\gamma'\|_{L^\infty(I_1)}$
  and $f_\star:=\|f'\|_{L^\infty(I_1)}$
  we obtain from Lemma \ref{lem3},
  (\ref{7.7}), (\ref{7.4}) and (\ref{7.2})
  that since $\na |\na v|^2 = 2D^2 v\cdot\na v$, and since $|\Del v| \le \sqrt{n} |D^2 v|$,
  \bea{7.76}
	& & \hs{-30mm}
	\frac{1}{p} \frac{d}{dt} \io |\na v|^p
	+ \gamma_\star \io |\na v|^{p-2} |D^2 v|^2 \nn\\
	&\le& - \frac{1}{2} \io \gamma'(\Theta) |\na v|^{p-2} \na\Theta\cdot\na |\na v|^2
	- \io \gamma'(\Theta) (\na\Theta\cdot\na v) \na v\cdot\na |\na v|^{p-2} \nn\\
	& & + 2a\io |\na v|^p
	+ a^{p+1} \io |\na u|^p \nn\\
	& & - \io \big( f'(\Theta)\cdot\na\Theta\big) \na v \cdot \na |\na v|^{p-2}
	- \io \big( f'(\Theta)\cdot\na\Theta\big) |\na v|^{p-2} \Del v  \nn\\
	&\le& \gamma_{\star\star} \io |\na\Theta| \cdot |\na v|^{p-1} |D^2 v|
	+ (p-2) \gamma_{\star\star} \io |\na \Theta| \cdot |\na v|^{p-1} |D^2 v| \nn\\
	& & + 2a\io |\na v|^p
	+ a^{p+1} \io |\na u|^p \nn\\
	& & + (p-2+\sqrt{n}) f_\star \io |\na\Theta| |\na v|^{p-2} |D^2 v| \nn\\
	&\le& \frac{\gamma_\star}{4} \io |\na v|^{p-2} |D^2 v|^2
	+  \frac{c_2 \gamma_{\star\star}^2}{\gamma_\star} \io |\na\Theta|^2 |\na v|^p \nn\\
	& & + 2a\io |\na v|^p
	+ a^{p+1} \io |\na u|^p  \nn\\
	&+& \frac{\gamma_\star}{4} \io |\na v|^{p-2} |D^2 v|^2
	+  \frac{c_3 f_\star^2}{\gamma_\star} \io |\na\Theta|^2 |\na v|^{p-2}
	\qquad \mbox{for all } t\in (0,T)
  \eea
  because of Young's inequality, where
  $c_2\equiv c_2(p):=(p-1)^2$ and $c_3\equiv c_3(p):=(p-2+\sqrt{n})^2$.
  When combined with (\ref{4.1}), this shows that
  \bea{7.777}
	& & \hs{-30mm}
	\frac{d}{dt} \bigg\{ \frac{1}{p} \io |\na v|^p
	+ \frac{8 a^{p-1} \gamma_\star \del_1}{p} \io |\na u|^p \bigg\}
	+ \frac{\gamma_\star}{2} \io |\na v|^{p-2} |D^2 v|^2
	+ 4a^p \gamma_\star \del_1 \io |\na u|^p \nn\\
	&\le& (2a+2^{p+2} \gamma_\star \del_1) \io |\na v|^p
	+ a^{p+1} \io |\na u|^p \nn\\
	& &
	+  \frac{c_2 \gamma_{\star\star}^2}{\gamma_\star} \io |\na\Theta|^2 |\na v|^p
	+  \frac{c_3 f_\star^2}{\gamma_\star} \io |\na\Theta|^2 |\na v|^{p-2}
	\qquad \mbox{for all } t\in (0,T),
  \eea
  where we note that thanks to (\ref{7.6}) and (\ref{7.5}), from (\ref{7.7}) it follows that
  \bea{7.77}
	(2a+2^{p+2}\del_1\gamma_\star) \io |\na v|^p
	\le 4(1+2^p)\del_1 \gamma_\star \io |\na v|^p
	\le \frac{\gamma_\star}{8c_1} \io |\na v|^p
	\le \frac{\gamma_\star}{8} \io |\na v|^{p-2} |D^2 v|^2
  \eea
  for all $t\in (0,T)$, as well as $a^{p+1} \le 2 a^p \gamma_\star \del_1$, that is,
  \be{7.78}
	a^{p+1} \io |\na u|^p \le 2a^p \gamma_\star \del_1 \io |\na u|^p
	\qquad \mbox{for all } t\in (0,T).
  \ee
  Moreover, Young's inequality ensures that
  \bea{7.79}
	\frac{c_2 \gamma_{\star\star}^2}{\gamma_\star}
	\io |\na\Theta|^2 |\na v|^p
	\le \io |\na v|^{p+2}
	+ c_2^\frac{p+2}{2} \gamma_{\star\star}^{p+2} \gamma_\star^{-\frac{p+2}{2}} \io |\na\Theta|^{p+2}
	\qquad \mbox{for all } t\in (0,T),
  \eea
  and that, again due to (\ref{7.5}),
  \bea{7.799}
	\frac{c_3 f_\star^2}{\gamma_\star} \io |\na\Theta|^2 |\na v|^{p-2}
	&=& \io \Big\{ \frac{\gamma_\star}{8 c_1} |\na v|^p \Big\}^\frac{p-2}{p} \cdot \Big(\frac{8 c_1}{\gamma_\star}\Big)^\frac{p-2}{p}
		\cdot \frac{c_3 f_\star^2}{\gamma\star} |\na\Theta|^2 \nn\\
	&\le& \frac{\gamma_\star}{8 c_1} \io |\na v|^p
	+ (8 c_1)^\frac{p-2}{2} c_3^\frac{p}{2} \gamma_\star^{1-p} f_\star^p \io |\na \Theta|^p \nn\\
	&\le& \frac{\gamma_\star}{8} \io |\na v|^{p-2} |D^2 v|^2
	+ c_4 \gamma_\star^{1-p} f_\star^p \io |\na \Theta|^p
	\qquad \mbox{for all } t\in (0,T)
  \eea
  with $c_4\equiv c_4(p):=(8 c_1)^\frac{p-2}{2} c_3^\frac{p}{2}$.
  Collecting (\ref{7.777})-(\ref{7.799}), we thus obtain that
  \bas
	& & \hs{-20mm}
	\frac{d}{dt} \bigg\{ \io |\na v|^p + 8 a^{p-1} \gamma_\star\del_1 \io |\na u|^p \bigg\}
	+ \frac{p\gamma_\star}{4} \io |\na v|^{p-2} |D^2 v|^2
	+ 2 p a^p \gamma_\star \del_1 \io |\na u|^p \nn\\
	&\le& p c_4 \gamma_\star^{1-p} f_\star^p \io |\na\Theta|^p
	+ p c_2^\frac{p+2}{2} \gamma_{\star\star}^{p+2} \gamma_\star^{-\frac{p+2}{2}} \io |\na \Theta|^{p+2}
	+ p \io |\na v|^{p+2}
	\qquad \mbox{for all } t\in (0,T),
  \eas
  and may hence conclude as intended, because once more by (\ref{7.5}),
  \bas
	\frac{p\gamma_\star}{4} \io |\na v|^{p-2} |D^2 v|^2
	\ge
	\frac{p\gamma_\star}{8} \io |\na v|^{p-2} |D^2 v|^2
	+ \frac{p\gamma_\star}{8 c_1} \io |\na v|^p
  \eas
  for all $t\in (0,T)$.
\qed
Parallel to Lemma \ref{lem7}, also the following statement on a similarly structured evolution property of
$t\mapsto \io |\na \Theta|^p$ relies on convexity of $\Om$ in dropping the last summand in (\ref{5.1}).
\begin{lem}\label{lem8}
  Let $n\ge 1$ and $\Om$ be convex, and let $p\ge 2$.
  Then there exist $k_2=k_2(p)>0$ and $K_2=K_2(p)>0$ such that if
  (\ref{gf}) and (\ref{init}) hold,
  and if $I\subset [0,\infty)$ is an interval such that
  \be{8.001}
	\Theta(x,t)\in I
	\qquad \mbox{for all $x\in\Om$ and } t\in (0,T)
  \ee
  with some $T\in (0,\tm]$, then
  \bea{8.1}
	& & \hs{-30mm}
	\frac{d}{dt}
	\io |\na\Theta|^p
	+ k_2 D \io |\na\Theta|^{p-2} |D^2\Theta|^2
	+ k_2 D \io |\na\Theta|^p \nn\\
	&\le& K_2 D^{1-p} \|F\|_{L^\infty(I)}^p \io |\na v|^p
	+ K_2 D^{1-p} a^p \|F\|_{L^\infty(I)}^p \io |\na u|^p  \nn\\
	& & + K_2 D^{-\frac{p+2}{4}} \|\Gamma\|_{L^\infty(I)}^\frac{p+2}{2} \io |\na v|^{p+2}
	+ K_2 D^{-\frac{p+2}{4}} \|\Gamma\|_{L^\infty(I)}^\frac{p+2}{2} a^{p+2} \io |\na u|^{p+2} \nn\\
	& & + 2 \io |\na\Theta|^{p+2}
	\qquad \mbox{for all } t\in (0,T).
  \eea
\end{lem}
\proof
  Once more by means of Lemma \ref{lem6},
  we fix $c_1=c_1(p)>0$ such that
  \be{8.02}
	\io |\na\vp|^p
	\le c_1 \io |\na\vp|^{p-2} |D^2\vp|^2
	\qquad \mbox{for all $\vp\in C^2(\bom)$ such that $\frac{\pa\vp}{\pa\nu}=0$ on $\pO$,}
  \ee
  and assume that $a,D,\gamma,\Gamma,f,F,I$ and $T$ have the listed properties.
  Again recalling \cite{lions_ARMA} in asserting that $\frac{\pa |\na\Theta|^2}{\pa\nu}\le 0$ on $\pO\times (0,T)$ by convexity
  of $\Om$, from Lemma \ref{lem5} and (\ref{8.02}) we readily obtain that if we let $K_0=K_0(p)$ be as provided by Lemma \ref{lem5}
  and write $c_2\equiv c_2(p):=\frac{p}{4}$, $c_3\equiv c_3(p):=\frac{p}{4c_1}$ as well as
  $\Gamma_\star:=\|\Gamma\|_{L^\infty(I)}$ and
  $F_\star:=\|F\|_{L^\infty(I)}$, we then have
  \bea{8.2}
	& & \hs{-30mm}
	\frac{d}{dt} \io |\na\Theta|^p
	+ c_2 D \io |\na\Theta|^{p-2} |D^2\Theta|^2
	+ c_3 D \io |\na\Theta|^p  \nn\\
	&\le& \frac{pK_0 \Gamma_\star^2}{D} \io |\na\Theta|^{p-2} |\na v|^4
	+ \frac{pK_0 a^4 \Gamma_\star^2}{D} \io |\na\Theta|^{p-2} |\na u|^4 \nn\\
	& & + \frac{pK_0 F_\star^2}{D} \io |\na\Theta|^{p-2} |\na v|^2
	+ \frac{pK_0 a^2 F_\star^2}{D} \io |\na\Theta|^{p-2} |\na u|^2
  \eea
  for all $t\in (0,T)$.
  Here since $p\ge 2$, Young's inequality says that
  \bas
	\frac{pK_0 \Gamma_\star^2}{D} \io |\na\Theta|^{p-2} |\na v|^4
	\le \io |\na\Theta|^{p+2}
	+ \Big(\frac{pK_0 \Gamma_\star^2}{D} \Big)^\frac{p+2}{4} \io |\na v|^{p+2}
  \eas
  as well as
  \bas
	\frac{pK_0 a^4 \Gamma_\star^2}{D} \io |\na\Theta|^{p-2} |\na u|^4
	\le \io |\na\Theta|^{p+2}
	+ \Big(\frac{pK_0 a^4 \Gamma_\star^2}{D} \Big)^\frac{p+2}{4} \io |\na u|^{p+2}
  \eas
  for all $t\in (0,T)$, and that
  \bas
	\frac{pK_0 F_\star^2}{D} \io |\na\Theta|^{p-2} |\na v|^2
	&=& \io \Big\{ \frac{c_4 D}{4} |\na\Theta|^p \Big\}^\frac{p-2}{p} \cdot \Big(\frac{4}{c_4 D}\Big)^\frac{p-2}{p}
		\cdot \frac{pK_0 F_\star^2}{D} |\na v|^2 \nn\\
	&\le& \frac{c_4 D}{4} \io |\na\Theta|^p
	+ \Big(\frac{4}{c_4}\Big)^\frac{p-2}{2} \cdot (pK_0)^\frac{p}{2} D^{1-p} F_\star^p \io |\na v|^p
  \eas
  and, similarly,
  \bas
	\frac{pK_0 a^2  F_\star^2}{D} \io |\na\Theta|^{p-2} |\na u|^2
	\le \frac{c_4 D}{4} \io |\na\Theta|^p
	+ \Big(\frac{4}{c_4}\Big)^\frac{p-2}{2} \cdot (pK_0)^\frac{p}{2} D^{1-p} a^p F_\star^p \io |\na u|^p
  \eas
  for all $t\in (0,T)$.
  Therefore, the claim results from (\ref{8.2}).
\qed
Now the superlinear contributions to the right-hand sides in (\ref{7.3}) and (\ref{8.1}) can suitably be controlled
on the basis of the following interpolation inequality.
\begin{lem}\label{lem99}
  Let $n\ge 1$ and $p\ge 2$ be such that $p>n$.
  Then there exists $K_3(p)>0$ with the property that whenever $\mu>0$,
  \be{99.1}
	\io |\na\vp|^{p+2} \le \mu \io |\na\vp|^{p-2} |D^2\vp|^2
	+ \mu \io |\na\vp|^p
	+ K_3(p) \mu^{-\frac{n}{p-n}} \cdot\bigg\{ \io |\na\vp|^p \bigg\}^\frac{p+2-n}{p-n}
	\qquad \mbox{for all } \vp\in C^2(\bom).
  \ee
\end{lem}
\proof
  As $p\ge n-2$, a Gagliardo-Nirenberg inequality yields $c_1=c_1(p)>0$ such that
  \bas
	\|\phi\|_{L^\frac{2(p+2)}{p}(\Om)}^\frac{2(p+2)}{p}
	\le c_1\|\na\phi\|_{L^2(\Om)}^\frac{2n}{p} \|\phi\|_{L^2(\Om)}^\frac{2(p+2-n)}{p}
	+ c_1 \|\phi\|_{L^2(\Om)}^\frac{2(p+2)}{p}
	\qquad \mbox{for all } \phi\in W^{1,2}(\Om),
  \eas
  which when applied to $\phi:=|\na\vp|^\frac{p}{2}$ implies that
  \bas
	\io |\na\vp|^{p+2}
	\le c_1\Big\| \na |\na\vp|^\frac{p}{2} \Big\|_{L^2(\Om)}^\frac{2n}{p}
		\Big\| |\na\vp|^\frac{p}{2} \Big\|_{L^2(\Om)}^\frac{2(p+2-n)}{p}
	+ c_1 \Big\| |\na\vp|^\frac{p}{2} \Big\|_{L^2(\Om)}^\frac{2(p+2)}{p}.
  \eas
  Using that actually $p>n$ and thus $\frac{n}{p}<1$, we may invoke Young's inequality to see that hence
  \bea{99.2}
	\io |\na\vp|^{p+2}
	&\le& \bigg\{ \frac{4\mu}{p^2} \Big\| \na |\na\vp|^\frac{p}{2} \Big\|_{L^2(\Om)}^2 \bigg\}^\frac{n}{p}
		\cdot \Big(\frac{p^2}{4\mu}\Big)^\frac{n}{p} c_1 \Big\| |\na\vp|^\frac{p}{2} \Big\|_{L^2(\Om)}^\frac{2(p+2-n)}{p}
	+ c_1 \Big\| |\na\vp|^\frac{p}{2} \Big\|_{L^2(\Om)}^\frac{2(p+2)}{p} \nn\\
	&\le& \frac{4\mu}{p^2} \Big\| \na |\na\vp|^\frac{p}{2} \Big\|_{L^2(\Om)}^2
	+ \Big(\frac{p^2}{4\mu}\Big)^\frac{n}{p-n} c_1^\frac{p}{p-n}
		\Big\| |\na\vp|^\frac{p}{2} \Big\|_{L^2(\Om)}^\frac{2(p+2-n)}{p-n}
	+ c_1 \Big\| |\na\vp|^\frac{p}{2} \Big\|_{L^2(\Om)}^\frac{2(p+2)}{p} \nn\\
	&=& \mu \io |\na\vp|^{p-2} |D^2\vp|^2
	+ c_2 \mu^{-\frac{n}{p-n}} \cdot \bigg\{ \io |\na\vp|^p \bigg\}^\frac{p+2-n}{p-n}
	+ c_1 \cdot\bigg\{ \io |\na\vp|^p \bigg\}^\frac{p+2}{p},
  \eea
  where $c_2\equiv c_2(p):=(\frac{p^2}{4})^\frac{n}{p-n} c_1^\frac{p}{p-n}$.
  Here the fact that $0<\frac{2}{p}<\frac{2}{p-n}$ enables us to once more utilize Young's inequality to verify that
  \bas
	c_1 \cdot\bigg\{ \io |\na\vp|^p \bigg\}^\frac{p+2}{p}
	&=& \bigg\{ \mu \io |\na\vp|^p \bigg\} \cdot \Bigg\{ 1\cdot\frac{c_1}{\mu} \cdot
		\bigg\{ \io |\na\vp|^p \bigg\}^\frac{2}{p} \Bigg\} \nn\\
	&\le& \bigg\{ \mu \io |\na\vp|^p \bigg\} \cdot \Bigg\{ 1 + \Big(\frac{c_1}{\mu})^\frac{p}{p-n} \cdot
		\bigg\{ \io |\na\vp|^p \bigg\}^\frac{2}{p-n} \Bigg\} \nn\\
	&=& \mu \io |\na v|^p
	+ c_1^\frac{p}{p-n} \mu^{-\frac{n}{p-n}} \cdot \bigg\{ \io |\na\vp|^p \bigg\}^\frac{p+2-n}{p-n}.
  \eas
  Consequently, (\ref{99.2}) implies (\ref{99.1}) with $K_3(p):=c_2+c_1^\frac{p}{p-n}$.
\qed
Indeed, by means of Lemma \ref{lem99} we can couple Lemma \ref{lem7} to Lemma \ref{lem8} and thereby derive
a superlinearly forced autonomous ODI
for a combination of expressions therein
in the style of (\ref{en1})-(\ref{energy}),
at least within settings of conveniently small data:
\begin{lem}\label{lem9}
  Let $n\ge 1$ and $\Om$ be convex, let $p\ge 2$ be such that $p>n$, and let $\del_1(p)\in (0,1)$
  be as provided by Lemma \ref{lem7}.
  Then there exists $\del(p)\in (0,\del_1(p)]$ such that if $a,D,\gamma,\Gamma,f$ and $F$ satisfy (\ref{gf}) and
  (\ref{7.1}) as well as (\ref{9.01}) with some $\Theta_\star\ge 0$,
  then taking $\eps_1=\eps_1(p,a,D,\gamma,\Gamma,f,F,\Theta_\star)$ as in Lemma \ref{lem7} one can find
  $\eps_2=\eps_2(p,a,D,\gamma,\Gamma,f,F,\Theta_\star)\in (0,\eps_1]$,
  $\eta_0=\eta_0(p,a,D,\gamma,\Gamma,f,F,\Theta_\star)>0$,
  $\kappa=\kappa(p,a,D,\gamma,\Gamma,f,F,\Theta_\star)>0$ and
  $K_4=K_4(p,a,D,\gamma,\Gamma,f,F,\Theta_\star)>0$ with the following property:
  If (\ref{init}) holds with
  \be{9.1}
	\io |\na v_0|^p
	+ \io |\na u_0|^p
	+ \io |\na u_0|^{p+2}
	+ \io |\na\Theta_0|^p
	\le \eta^p
  \ee
  for some $\eta\in (0,\eta_0)$, and if $T\in (0,\tm]$ is such that (\ref{7.2}) holds with some $\eps\in (0,\eps_2]$, then
  \be{9.2}
	\io |\na v(\cdot,t)|^p
	+ \io |\na u(\cdot,t)|^p
	+ \io |\na u(\cdot,t)|^{p+2}
	+ \io |\na \Theta(\cdot,t)|^p
	\le K_4 \eta^p e^{-\kappa t}
	\qquad \mbox{for all } t\in (0,T)
  \ee
  and
  \be{9.3}
	\int_0^t \io |\na v|^{p-2} |D^2 v|^2
	+ \int_0^t \io |\na \Theta|^{p-2} |D^2 \Theta|^2
	\le K_4 \eta^p
	\qquad \mbox{for all } t\in (0,T).
  \ee
\end{lem}
\proof
  We fix $k_i=k_i(p,a,D,\gamma,\Gamma,f,F,\Theta_\star)$ and $K_i=K_i(p,a,D,\gamma,\Gamma,f,F,\Theta_\star)$, $i\in\{1,2\}$,
  as in Lemma \ref{lem7} and Lemma \ref{lem8}, and set
  \bas
	\del\equiv \del(p):=\min \bigg\{ \del_1(p) \, , \, \Big\{\frac{k_1 k_2}{8 K_1 K_2}\Big\}^\frac{1}{p} \bigg\}.
  \eas
  Then whenever (\ref{7.1}) and (\ref{9.01}) hold with some $\Theta_\star\ge 0$,
  by continuity of $f', F$ and $\gamma$ we can fix $\eps_2=\eps_2(p,a,D,\gamma,\Gamma,F,f,\Theta_\star)\in (0,\eps_1]$
  in such a way that for
  $I:=\{\xi\ge 0 \ | \ |\xi-\Theta_\star|\le \eps_2\}$,
  $f_\star:=\|f'\|_{L^\infty(I)}$,
  $F_\star:=\|F\|_{L^\infty(I)}$ and
  $\gamma_\star:=\inf_{\xi\in I} \gamma(\xi)$ we have
  \bas
	\frac{f_\star^p F_\star^p}{D^p \gamma_\star^p} < \frac{k_1 k_2}{4K_1 K_2}.
  \eas
  As thus
  \bas
	\frac{K_1 \gamma_\star^{1-p} f_\star^p}{\frac{k_2 D}{2}} <
	\frac{\frac{k_1 \gamma_\star}{2}}{K_2 D^{1-p} F_\star^p},
  \eas
  it is possible to choose $A=A(p,a,D,\gamma,\Gamma,F,f,\Theta_\star)>0$ such that
  \be{9.11}
	K_1 \gamma_\star^{1-p} f_\star^p \le A\cdot\frac{k_2 D}{2},
  \ee
  and that
  \be{9.12}
	A\cdot K_2 D^{1-p} F_\star^p \le \frac{k_1\gamma_\star}{2}
  \ee
  as well as, equivalently,
  \be{9.13}
	A\cdot K_2 a^p D^{1-p} F_\star^p \le \frac{k_1 a^p \gamma_\star}{2}.
  \ee
  We thereupon let $B=B(p,a,D,\gamma,\Gamma,F,f,\Theta_\star)>0$ be large enough fulfilling
  \be{9.14}
	A\cdot K_2 D^{-\frac{p+2}{4}} a^{p+2} \Gamma_\star^\frac{p+2}{2}
	\le B\cdot\frac{(p+2) a}{4}
  \ee
  with $\Gamma_\star:=\|\Gamma\|_{L^\infty(I)}$, and taking $K_3=K_3(p)$ from Lemma \ref{lem99} we define
  \be{9.15}
	\kappa\equiv\kappa(p,a,D,\gamma,\Gamma,F,f,\Theta_\star):=\min \bigg\{ \frac{k_1\gamma_\star}{8} \, , \, \frac{k_1 a}{32\del_1(p)} \, , \, \frac{k_2 D}{8}
		\, , \, \frac{(p+2)a}{8} \bigg\}
  \ee
  and $\lam\equiv \lam(p):=\frac{p+2-n}{p-n}$
  as well as the numbers $c_i\equiv c_i(p,a,D,\gamma,\Gamma,F,f,\Theta_\star)$, $i\in\{1,...,6\}$, according to
  \be{9.16}
	c_1:=K_1+A\cdot K_2 D^{-\frac{p+2}{4}} \Gamma_\star^\frac{p+2}{2} + B\cdot 2^{p+1} (p+2) a^{-p-1}
	\qquad \mbox{and} \qquad
	c_2:=K_1 \gamma_\star^{-\frac{p+2}{2}} \gamma_{\star\star}^{p+2} + 2A,
  \ee
  to
  \be{9.17}
	c_3:=c_1K_3\cdot \Big(\frac{4c_1}{k_1\gamma_\star} \Big)^\frac{n}{p-n}
	\qquad \mbox{and} \qquad
	c_4:=c_2 K_3\cdot \Big(\frac{4c_2}{A\cdot k_2 D} \Big)^\frac{n}{p-n},
  \ee
  and to
  \be{9.177}
	c_5:=c_3+\frac{c_4}{A^\lam}
	\qquad \mbox{and} \qquad
	c_6:=\max \big\{ 1 \, , \, 8a^{p-1} \gamma_\star \del_1(p) \, , \, A \, , \, B \big\}.
  \ee
  We finally introduce
  \be{9.18}
	\eta_0\equiv \eta_0(p,a,D,\gamma,\Gamma,F,f,\Theta_\star):=\Big(\frac{\kappa}{c_5 c_6^{\lam-1}}\Big)^\frac{1}{p(\lam-1)},
  \ee
  and from now on we assume that (\ref{init}) and (\ref{9.1}) hold with some $\eta\in (0,\eta_0]$, and that $T\in (0,\tm]$
  is such that (\ref{7.2}) is satisfied with some $\eps\in (0,\eps_2)$, we let
  \be{9.19}
	y(t):=
	\io |\na v(\cdot,t)|^p
	+ 8a^{p-1} \gamma_\star \del_1(p) \io |\na u(\cdot,t)|^p
	+ A \io |\na\Theta(\cdot,t)|^p
	+ B \io |\na u(\cdot,t)|^{p+2},
	\quad t\in [0,T),
  \ee
  as well as
  \be{9.20}
	h(t):=\frac{k_1\gamma_\star}{2} \io |\na v(\cdot,t)|^{p-2} |D^2 v(\cdot,t)|^2
	+ \frac{A\cdot k_2 D}{2} \io |\na\Theta(\cdot,t)|^{p-2}|D^2\Theta(\cdot,t)|^2,
	\qquad t\in (0,T),
  \ee
  and combine (\ref{7.3}) with (\ref{8.1}) and (\ref{4.2}) to see that
  \bas
	& & \hs{-12mm}
	y'(t) + h(t)
	+ \frac{k_1\gamma_\star}{2} \io |\na v|^{p-2} |D^2 v|^2
	+ \frac{A\cdot k_2 D}{2} \io |\na\Theta|^{p-2} |D^2\Theta|^2 \\
	& & + k_1 \gamma_\star \io |\na v|^p
	+ k_1 a^p \gamma_\star \io |\na u|^p \\
	& & + A\cdot k_2 D \io |\na\Theta|^p
	+ B\cdot \frac{(p+2)a}{2} \io |\na u|^{p+2} \nn\\
	&\le& \bigg\{ K_1 \gamma_\star^{1-p} f_\star^p \io |\na\Theta|^p
	+ K_1 \gamma_\star^{-\frac{p+2}{2}} \gamma_{\star\star}^{p+2} \io |\na\Theta|^{p+2} + K_1 \io |\na v|^{p+2} \bigg\} \nn\\
	& & + \bigg\{ A\cdot K_2 D^{1-p} F_\star^p \io |\na v|^p
	+ A\cdot K_2 D^{1-p} a^p F_\star^p \io |\na u|^p \nn\\
	& & \hs{6mm}
	+ A\cdot K_2 D^{-\frac{p+2}{4}} \Gamma_\star^\frac{p+2}{2} \io |\na v|^{p+2}
	+ A\cdot K_2 D^{-\frac{p+2}{4}} a^{p+2} \Gamma_\star^\frac{p+2}{2} \io |\na u|^{p+2}
	+ 2A \io |\na\Theta|^{p+2} \bigg\} \nn\\
	& & + B\cdot 2^{p+1} (p+2) a^{-p-1} \io |\na v|^{p+2}
	\qquad \mbox{for all } t\in (0,T).
  \eas
  In view of (\ref{9.11})-(\ref{9.16}), this readily implies that
  \bas
	& & \hs{-20mm}
	y'(t) + h(t) + 2\kappa y(t)
	+ \frac{k_1 \gamma_\star}{4} \io |\na v|^p
	+ \frac{A\cdot k_2 D}{4} \io |\na\Theta|^p \nn\\
	&\le& c_1 \io |\na v|^{p+2} + c_2 \io |\na\Theta|^{p+2}
	\qquad \mbox{for all } t\in (0,T),
  \eas
  and here two applications of Lemma \ref{lem99} show that
  \bas
	c_1 \io |\na v|^{p+2}
	\le \frac{k_1 \gamma_\star}{4} \io |\na v|^{p-2} |D^2 v|^2
	+ \frac{k_1 \gamma_\star}{4} \io |\na v|^p
	+ c_1 K_3 \cdot \Big(\frac{4c_1}{k_1\gamma_\star}\Big)^\frac{n}{p-n} \cdot \bigg\{ \io |\na v|^p \bigg\}^\lam
  \eas
  and
  \bas
	c_2 \io |\na \Theta|^{p+2}
	\le \frac{A\cdot k_2 D}{4} \io |\na\Theta|^{p-2} |D^2\Theta|^2
	+ \frac{A\cdot k_2 D}{4} \io |\na\Theta|^p
	+ c_2 K_3 \cdot \Big(\frac{4c_2}{A\cdot k_2 D}\Big)^\frac{n}{p-n} \cdot \bigg\{ \io |\na \Theta|^p \bigg\}^\lam
  \eas
  for all $t\in (0,T)$.
  In line with (\ref{9.17}), (\ref{9.177}) and (\ref{9.19}), we thus infer that
  \bea{9.21}
	y'(t) + 2\kappa y(t) + \frac{1}{2} h(t)
	&\le& c_3 \cdot \bigg\{ \io |\na v|^p \bigg\}^\lam
	+ c_4\cdot\bigg\{ \io |\na\Theta|^p \bigg\}^\lam \nn\\
	&\le& c_3 y^\lam(t)
	+ c_4 \cdot \Big(\frac{y(t)}{A}\Big)^\lam
	= c_5 y^\lam(t)
	\qquad \mbox{for all } t\in (0,T),
  \eea
  while, on the other hand, for
  \bas
	\oy(t):=c_6 \eta^p e^{-\kappa t},
	\qquad t\ge 0,
  \eas
  we have
  \bas
	\oy'(t) + 2\kappa \oy(t) + \frac{1}{2} h(t) - c_5 \oy^\lam(t)
	&\ge& \oy'(t) + 2\kappa \oy(t) - c_5 \oy^\lam(t) \nn\\
	&=& - c_6\kappa \eta^p e^{-\kappa t}
	+ 2c_6 \kappa \eta^p e^{-\kappa t}
	- c_5 c_6^\lam \eta^{p\lam} e^{-\kappa\lam t} \\
	&=& c_6\kappa \eta^p e^{-\kappa t} \cdot \bigg\{ 1 - \frac{c_5 c_6^{\lam-1}}{\kappa} \eta^{p(\lam-1)} e^{-\kappa(\lam-1)t}
		\bigg\} \\
	&\ge& c_6 \kappa \eta^p e^{-\kappa t} \cdot \bigg\{ 1 - \frac{c_5 c_6^{\lam-1}}{\kappa} \eta_0^{p(\lam-1)} \bigg\}
	\ge 0
	\qquad \mbox{for all } t>0
  \eas
  due to (\ref{9.18}).
  As furthermore from (\ref{9.1}) and (\ref{9.177}) we know that
  \bas
	y(0) \le c_6 \cdot \bigg\{ \io |\na v_0|^p + \io |\na u_0|^p + \io |\na\Theta_0|^p + \io |\na u_0|^{p+2} \bigg\}
	\le c_6\eta^p =\oy(0),
  \eas
  an ODE comparison argument shows that
  \be{9.22}
	y(t)\le \oy(t)
	\qquad \mbox{for all } t\in (0,T).
  \ee
  Thereupon, an integration in (\ref{9.21}) reveals that
  \bas
	\frac{1}{2} \int_0^t h(s) ds
	&\le& y(0) + c_5 \int_0^t y^\lam(s) ds \\
	&\le& c_6 \eta^p + c_5 c_6^\lam \eta^{p\lam} \int_0^\infty e^{-\lam\kappa s} ds
	\qquad \mbox{for all } t\in (0,T),
  \eas
  which together with (\ref{9.22}) yields both (\ref{9.2}) and (\ref{9.3}) with some adequately large positive
  $K_4=K_4(p,a,D,\gamma,\Gamma,f,F,\Theta_\star)$.
\qed
As a final preparation, let us record an elementary consequence of the elementary conservation feature
$\frac{d^2}{dt^2} \io u=0$ for solutions of (\ref{0}), here translated to the framework of (\ref{0v}).
\begin{lem}\label{lem_mass}
  Let $n\ge 1$ and $\Omega$ be an arbitrary bounded domain with smooth boundary, and suppose that (\ref{gf}) is satisfied.
  Then assuming (\ref{init}), we have
  \be{mass}
	\io \big(v(\cdot,t)-au(\cdot,t)\big)=\io (v_0-au_0)
	+ \int_0^t \int_{\pO} f(\Theta)\cdot\nu
	\qquad \mbox{for all } t\in (0,\tm);
  \ee
  in particular,
  for all $p\ge 1$ there exists $C=C(p)>0$ such that
  whenever (\ref{init}) holds, it follows that
  \be{m1}
	\|v(\cdot,t) -a u(\cdot,t)\|_{W^{1,p}(\Om)}
	\le C \|\na v(\cdot,t)-a \na u(\cdot,t)\|_{L^p(\Om)}
	+ C \cdot \bigg| \io (v_0-au_0)\bigg|
	+ C \int_0^t \|f(\Theta(\cdot,s))\|_{L^\infty(\Om)} ds
  \ee
  for all $t\in (0,\tm)$.
\end{lem}
\proof
  From (\ref{0v}) we immediately obtain that
  \bas
	\frac{d}{dt} \io v - a\frac{d}{dt} \io u
	= \io \na\cdot f(\Theta)
	= \int_{\pO} f(\Theta)\cdot\nu
	\qquad \mbox{for all } t\in (0,\tm).
  \eas
  This implies (\ref{mass}), whereupon (\ref{m1}) readily results by means of a Poincar\'e-type inequality, according to
  which, namely, there exists $c_1=c_1(p)>0$ such that
  \bas
	\bigg\| \vp - \frac{1}{|\Om|} \io \vp \bigg\|_{W^{1,p}(\Om)} \le c_1 \|\na \vp\|_{L^p(\Om)}
  \eas
  for all $\vp\in W^{1,p}(\Om)$.
\qed
In fact, this enables us to close our loop of arguments, and to particularly confirm that for small-data solutions
the second alternative in (\ref{ext}) cannot hold:
\begin{lem}\label{lem11}
  Let $n\ge 1$ and $\Om$ be convex, let $p\ge 2$ be such that $p>n$, and suppose that (\ref{gf}) as well as (\ref{7.1})
  and (\ref{9.01}) hold with some $\Theta_\star\ge 0$.
  Then taking $\eps_2=\eps_2(p,a,D,\gamma,\Gamma,f,F,\Theta_\star)$ and $\eta_0=\eta_0(p,a,D,\gamma,\Gamma,f,F,\Theta_\star)$
  from Lemma \ref{lem9}, given any $\eps\in (0,\frac{\eps_2}{2})$ one can find
  $\eta=\eta(\eps,p,a,D,\gamma,\Gamma,$ $f,F,\Theta_\star)\in (0,\eta_0)$ such that
  whenever (\ref{init}) and (\ref{9.1}) as well as
  \be{11.1}
	\|\Theta_0-\Theta_\star\|_{L^\infty(\Om)} \le \eps
  \ee
  hold, it follows that $\tm=\infty$, and that
  \be{11.2}
	\|\Theta(\cdot,t)-\Theta_\star\|_{L^\infty(\Om)} \le 2\eps
	\qquad \mbox{for all } t>0.
  \ee
\end{lem}
\proof
  According to known smoothing properties of the Neumann heat semigroup $(e^{t\Del})_{t\ge 0}$ on $\Om$, there exist
  $c_1=c_1(p,D)>0$ and $c_2=c_2(p,D)>0$ such that for all $t>0$ and any $\vp\in C^0(\bom)$ we have
  \be{11.3}
	\|e^{tD\Del}\vp\|_{L^\infty(\Om)}
	\le c_1 \big(1+t^{-\frac{n}{p}}\big) \|\vp\|_{L^\frac{p}{2}(\Om)}
  \ee
  and
  \be{11.33}
	\|e^{tD\Del}\vp\|_{L^\infty(\Om)}
	\le c_2 \big(1+t^{-\frac{n}{2p}}\big) \|\vp\|_{L^p(\Om)}.
  \ee
  Taking $K_4=K_4(p,a,D,\gamma,\Gamma,f,F,\Theta_\star)$ and $\kappa=\kappa(p,a,D,\gamma,\Gamma,f,F,\Theta_\star)$ from Lemma \ref{lem9},
  abbreviating $I:=\{\xi\ge 0 \ | \ |\xi-\Theta_\star| \le \eps_2\}$,
  $\Gamma_\star:=\|\Gamma\|_{L^\infty(I)}$ and
  $F_\star:=\|F\|_{L^\infty(I)}$,
  and noting that both
  \be{11.4}
	c_3\equiv c_3(p,a,D,\gamma,\Gamma,f,F,\Theta_\star):=\sup_{t>0} \int_0^t \big( 1+(t-s)^{-\frac{n}{p}}\big) e^{-\frac{2\kappa s}{p}} ds
  \ee
  and
  \be{11.44}
	c_4\equiv c_4(p,a,D,\gamma,\Gamma,f,F,\Theta_\star):=\sup_{t>0} \int_0^t \big( 1+(t-s)^{-\frac{n}{2p}}\big) e^{-\frac{\kappa s}{p}} ds
  \ee
  are finite due to the fact that $p>n>\frac{n}{2}$, given any $\eps\in (0,\frac{\eps_2}{2})$
  we choose $\eta=\eta(\eps,p,a,D,\gamma,\Gamma,f,F,\Theta_\star)\in (0,\eta_0)$ in such a way that
  \be{11.5}
	2c_1 c_3 \Gamma_\star (1+a^2) K_4^\frac{2}{p} \eta^2 \le \frac{\eps}{4}
	\quad \mbox{and} \quad
	c_2 c_4 F_\star (1+a) K_p^\frac{1}{p} \eta \le \frac{\eps}{4}.
  \ee
  We now suppose that (\ref{init}), (\ref{9.1}) and (\ref{11.1}) hold,
  and claim that the above selections ensure that
  \be{11.55}
	T:=\sup \Big\{ \wh{T}\in (0,\tm) \ \Big| \ \|\Theta(\cdot,t)-\Theta_\star\|_{L^\infty(\Om)} < 2\eps \mbox{ for all }
		t\in (0,\wh{T}) \Big\},
  \ee
  being a well-defined element of $(0,\tm] \subset (0,\infty]$ due to (\ref{11.1}), actually satisfies
  \be{11.6}
	T=\tm.
  \ee
  To see this, we observe that according to a Duhamel representation associated with the identity
  \bas
	(\Theta-\Theta_\star)_t = D\Del (\Theta-\Theta_\star) + \Gamma(\Theta) |\na v-a\na u|^2 + F(\Theta)(\na v-a\na u),
  \eas
  and by the maximum principle, (\ref{11.3}) and (\ref{11.33}), we have
  \bea{11.7}
	& & \hs{-20mm}
	\|\Theta(\cdot,t)-\Theta_\star\|_{L^\infty(\Om)} \nn\\
	&=& \bigg\| e^{tD\Del} (\Theta_0-\Theta_\star) + \int_0^t e^{(t-s)D\Del}
		\Big\{ \Gamma(\Theta(\cdot,s)) \big| \na v(\cdot,s)-a\na u(\cdot,s)\big|^2 \Big\} ds \nn\\
	& & \hs{40mm}
	+ \int_0^t e^{(t-s)D\Del} \Big\{ F(\Theta(\cdot,s)) \big(\na v(\cdot,s)-a\na u(\cdot,s)\big) \Big\} ds
		\bigg\|_{L^\infty(\Om)} \nn\\
	&\le& \| e^{tD\Del} (\Theta_0-\Theta_\star)\|_{L^\infty(\Om)}
	+ \int_0^t \Big\| e^{(t-s)D\Del}
		\Big\{ \Gamma(\Theta(\cdot,s)) \big| \na v(\cdot,s)-a\na u(\cdot,s)\big|^2 \Big\} ds \Big\|_{L^\infty(\Om)} \nn\\
	& & + \int_0^t \Big\| e^{(t-s)D\Del}
		\Big\{ F(\Theta(\cdot,s)) \big( \na v(\cdot,s)-a\na u(\cdot,s)\big) \Big\} ds \Big\|_{L^\infty(\Om)} \nn\\
	&\le& \|\Theta_0-\Theta_\star\|_{L^\infty(\Om)}
	+ c_1 \int_0^t \big(1+(t-s)^{-\frac{n}{p}}\big)
		\Big\| \Gamma(\Theta(\cdot,s)) \big| \na v(\cdot,s)-a\na u(\cdot,s)\big|^2 \Big\|_{L^\frac{p}{2}(\Om)} ds \nn\\
	& & + c_2 \int_0^t \big(1+(t-s)^{-\frac{n}{2p}}\big)
		\Big\| F(\Theta(\cdot,s)) \big( \na v(\cdot,s)-a\na u(\cdot,s)\big) \Big\|_{L^p(\Om)} ds \nn\\
	&\le& \|\Theta_0-\Theta_\star\|_{L^\infty(\Om)}
	+ c_1 \Gamma_\star \int_0^t \big(1+(t-s)^{-\frac{n}{p}}\big)
		\Big\| \big| \na v(\cdot,s)-a\na u(\cdot,s)\big|^2 \Big\|_{L^\frac{p}{2}(\Om)} ds \nn\\
	& & + c_2 F_\star \int_0^t \big(1+(t-s)^{-\frac{n}{2p}}\big)
		\big\| \na v(\cdot,s)-a\na u(\cdot,s)\big\|_{L^p(\Om)} ds
  \eea
  for all $t\in (0,T)$,
  because $\Theta(x,s)\in I$ for all $(x,s)\in \Om\times (0,T)$ by (\ref{11.55}) and the inequality $2\eps\le \eps_2$.
  Here since $\eta<\eta_0$, Lemma \ref{lem9} applies so as to ensure that
  \bas
	\Big\| |\na v-a\na u|^2 \Big\|_{L^\frac{p}{2}(\Om)}
	&=& \|\na v-a\na u\|_{L^p(\Om)}^2 \\
	&\le& 2\|\na v\|_{L^p(\Om)}^2
	+ 2a^2 \|\na u\|_{L^p(\Om)}^2 \\
	&\le& 2\cdot \big(K_4\eta^p e^{-\kappa t}\big)^\frac{2}{p}
	+ 2a^2 \cdot (K_4\eta^p e^{-\kappa t}\big)^\frac{2}{p} \\
	&=& 2(1+a^2) K_4^\frac{2}{p} \eta^2 e^{-\frac{2\kappa t}{p}}
  \eas
  and
  \bas
	\big\| \na v-a\na u\|_{L^p(\Om)}
	\le \|\na v\|_{L^p(\Om)}
	+ a\|\na u\|_{L^p(\Om)}
	\le (1+a) \cdot (K_4\eta^p e^{-\kappa t}\big)^\frac{1}{p}
	= (1+a)K_4^\frac{1}{p} \eta e^{-\frac{\kappa t}{p}}
  \eas
  for all $t\in (0,T)$,
  so that in line with (\ref{11.4}), (\ref{11.44}) and (\ref{11.5}),
  \bas
	& & \hs{-20mm}
	c_1 \Gamma_\star \int_0^t \big(1+(t-s)^{-\frac{n}{p}}\big)
		\Big\| \big| \na v(\cdot,s)-a\na u(\cdot,s)\big|^2 \Big\|_{L^\frac{p}{2}(\Om)} ds \\
	&\le& 2c_1 \Gamma_\star (1+a^2) K_4^\frac{2}{p} \eta^2
		\int_0^t \big(1+(t-s)^{-\frac{n}{p}}\big) e^{-\frac{2\kappa s}{p}} ds \\
	&\le& 2c_1 \Gamma_\star (1+a^2) K_4^\frac{2}{p} \eta^2 c_3 \\
	&\le& \frac{\eps}{4}
	\qquad \mbox{for all } t\in (0,T)
  \eas
  and
  \bas
	& & \hs{-20mm}
	c_2 F_\star \int_0^t \big(1+(t-s)^{-\frac{n}{2p}}\big)
		\big\| \na v(\cdot,s)-a\na u(\cdot,s)\big\|_{L^\frac{p}{2}(\Om)} ds \\
	&\le& c_2 F_\star (1+a) K_4^\frac{1}{p} \eta
		\int_0^t \big(1+(t-s)^{-\frac{n}{2p}}\big) e^{-\frac{\kappa s}{p}} ds \\
	&\le& c_2 F_\star (1+a) K_p^\frac{1}{p} \eta c_4 \\
	&\le& \frac{\eps}{4}
	\qquad \mbox{for all } t\in (0,T).
  \eas
  Thanks to (\ref{11.1}), from (\ref{11.7}) we thus infer that
  \bas
	\|\Theta(\cdot,t)-\Theta_\star\|_{L^\infty(\Om)} \le \eps+\frac{\eps}{4}+ \frac{\eps}{4}
	\qquad \mbox{for all } t\in (0,T),
  \eas
  which by continuity of $\Theta$ is compatible with (\ref{11.55}) only if indeed (\ref{11.6}) holds.\abs
  Thus, $|\Theta-\Theta_\star|\le 2\eps$ throughout $\Om\times (0,\tm)$, and to finally make sure that $\tm=\infty$,
  we only need to observe that if $\tm$ was finite, then
  (\ref{11.6}) together with Lemma \ref{lem9}
  and Lemma \ref{lem_mass}
  would ensure that also $\sup_{t\in (0,\tm)} \|v(\cdot,t)-au(\cdot,t)\|_{W^{1,p}(\Om)}<\infty$.
  Combined with the $L^\infty$ bound for $\Theta$ just established,
  this would contradict (\ref{ext}).
\qed
As a consequence of this and of our local theory for (\ref{0v}), any of the solutions addressed in Lemma \ref{lem9}
must be global in time and have the announced decay properties:\abs
\proofc of Theorem \ref{theo12}. \quad
  This readily results on applying Lemma \ref{lem_loc}, Lemma \ref{lem9} and Lemma \ref{lem11}, again to
  $(v_0,u_0,\Theta_0):=(u_{0t}+au_0,u_0,\Theta_0)$, noting that if (\ref{12.2}) holds, then
  \bas
	& & \hs{-30mm}
	\io |\na v_0|^p + \io |\na u_0|^p + \io |\na u_0|^{p+2} + \io |\na\Theta_0|^p \\
	&=& \io |\na u_{0t}+a\na u_0|^p + \io |\na u_0|^p + \io |\na u_0|^{p+2} + \io |\na\Theta_0|^p \\
	&\le& 2^{p-1} \io |\na u_{0t}|^p + (1+2^{p-1} a^p) \io |\na u_0|^p + \io |\na u_0|^{p+2} + \io |\na\Theta_0|^p \\
	&\le& \max \big\{ 2^{p-1} \, , \, 1+2^{p-1} a^p \big\} \cdot \eta^p
  \eas
  due to Young's inequality.
\qed

\bigskip

{\bf Acknowledgment.} \quad
The authors acknowledge support of the Deutsche Forschungsgemeinschaft (Project No. 444955436).
They moreover declare that they have no conflict of interest.\abs
{\bf Data availability statement.} \quad
Data sharing is not applicable to this article as no datasets were
generated or analyzed during the current study.

\end{document}